\documentclass[12pt, colorinlistoftodos]{amsart}
\usepackage{a4wide}
\usepackage{amssymb}
\usepackage{amsthm}
\usepackage{amsmath}
\usepackage{amscd}
\usepackage{verbatim}
\usepackage{bm}
\usepackage[mathscr]{eucal}
\usepackage[all]{xy}
\usepackage{todonotes}
\usepackage{enumitem}
\usepackage[pagebackref=true]{hyperref}
\usepackage{dsfont, stmaryrd}
\usepackage{mathrsfs}
\usepackage{tikz, tikz-cd, caption, tabu}
\usepackage{tkz-euclide, mathtools}
\numberwithin{equation}{section}

\theoremstyle{plain}
\newtheorem{theorem}{Theorem}[section]

\newtheorem{lemma}[theorem]{Lemma}
\newtheorem{proposition}[theorem]{Proposition}

\theoremstyle{definition}

\newtheorem{remark}[theorem]{Remark}
\newtheorem{example}[theorem]{Example}

\theoremstyle{remark}

\newcommand{\pmat}[4]{\begin{pmatrix}
                 #1 & #2\\
                 #3 & #4
\end{pmatrix}}
\newcommand{\smat}[4]{\left(\begin{smallmatrix}
                 #1 & #2\\
                 #3 & #4
\end{smallmatrix}\right)}

\newcommand{\lp}{\left (}
\newcommand{\rp}{\right )}

\newcommand{\Fc}{{\mathcal{F}}}

\newcommand{\Oc}{{\mathcal{O}}}
\newcommand{\Sc}{{\mathcal{S}}}
\newcommand{\Zb}{\mathbb{Z}}

\newcommand{\Qb}{\mathbb{Q}}
\newcommand{\SO}{{\mathrm{SO}}}
\newcommand{\SL}{{\mathrm{SL}}}
\newcommand{\GL}{{\mathrm{GL}}}

\newcommand{\af}{\mathfrak{a}}

\newcommand{\ef}{\mathfrak{e}}

\newcommand{\wf}{\mathfrak{w}}

\newcommand{\vf}{\mathfrak{v}}

\newcommand{\Nm}{{\mathrm{Nm}}}
\newcommand{\Ab}{\mathbb{A}}
\newcommand{\Db}{\mathbb{D}}
\newcommand{\Nb}{\mathbb{N}}
\newcommand{\Hb}{\mathbb{H}}
\newcommand{\Rb}{\mathbb{R}}
\newcommand{\Cb}{\mathbb{C}}

\newcommand{\Gm}{\mathbb{G}_m}
\newcommand{\ebf}{{\mathbf{e}}}

\newcommand{\tf}{\tilde{f}}

\newcommand{\slf}{{\mathfrak{sl}}}
\newcommand{\spf}{{\mathfrak{sp}}}

\newcommand{\zf}{\mathfrak{z}}

\newcommand{\reg}{\mathrm{reg}}

\newcommand{\Pc}{\mathcal{P}}
\newcommand{\Qc}{\mathcal{Q}}
\newcommand{\Rc}{\mathcal{R}}

\newcommand{\CT}{\mathrm{CT}}

\newcommand{\Mp}{\mathrm{Mp}}

\newcommand{\cha}{\mathrm{char}}

\newcommand{\GSpin}{\mathrm{GSpin}}
\newcommand{\Lz}{L_z}
\newcommand{\Rz}{R_z}
\newcommand{\He}{\mathrm{H}}

\renewcommand{\k}{E}

\newcommand{\BS}{\mathrm{BS}}
\newcommand{\BZ}{\mathrm{BZ}}
\renewcommand{\L}{\mathbb{L}}
\AtBeginDocument{%
   \def\MR#1{}
}

\begin{document}
\title{Laurent Expansions of Meromorphic Modular Forms}
\author[G.~ Bogo, Y.~Li, M.~Schwagenscheidt]{Gabriele Bogo, Yingkun Li, Markus Schwagenscheidt}
\address{Fachbereich Mathematik,
Technische Universit\"at Darmstadt, Schlossgartenstrasse 7, D--64289
Darmstadt, Germany}
\email{bogo@mathematik.tu-darmstadt.de}
\email{li@mathematik.tu-darmstadt.de}
\address{ETH Z\"urich Mathematics Dept., R\"amistrasse 101, CH-8092 Z\"urich, Switzerland}
\email{mschwagen@ethz.ch}

\subjclass[2020]{}
\thanks{
}

\begin{abstract}
  In this paper, we study the Laurent coefficients of meromorphic modular forms at CM points by giving two approaches of computing them.
  The first is a generalization of the method of Rodriguez--Villegas and Zagier, which expresses the Laurent coefficients as constant terms of a family of polynomials obtained through recursion.
  The second applies to meromorphic modular forms that are regularized theta lifts, and expresses their Laurent coefficients in terms of Fourier coefficients of harmonic Maass forms. 
\end{abstract}

\date{\today}
\maketitle

 \makeatletter
 \providecommand\@dotsep{5}
 \def\listtodoname{List of Todos}
 \def\listoftodos{\@starttoc{tdo}\listtodoname}
 \makeatother


\section{Introduction}

Many striking applications of modular forms rely on the fact that their Fourier coefficients are typically given by arithmetically interesting quantities.  Classical examples include
generalized divisors sums, representation numbers of quadratic forms, or the partition function. As a result, the Fourier coefficients of modular forms have been studied extensively since the 19th century and remain a central topic of research until today.

In contrast, the Laurent expansions of modular forms around points in the upper half-plane have received less attention. Here, the Laurent expansion of a holomorphic modular form $g$ of weight $k \in \Zb$ for $\SL_2(\Zb)$ at a point $z_0 \in \Hb$ is defined by the expansion
\[
  (1-w)^{-k}g(z) = \sum_{m \geq 0}
  a_m w^m,
\]
with the local variable $w = \frac{z-z_0}{z-\overline{z}_0}$ around $z_0$. Introducing the local variable $w$ has the advantage over the ``usual'' Laurent expansion (involving just powers of $z-z_0$) that the above expansion converges for all $z \in \Hb$. Moreover, the Laurent coefficients $a_m$ can be obtained as special values at $z = z_0$ of ``modular'' derivatives of $g$. If we define the weight $k$ raising operator $R_{k} = 2i\frac{\partial}{\partial z}+ky^{-1}$ and its iterated version
\[
R_{k}^m := R_{k+2m-2}\circ\cdots \circ R_{k+2}\circ R_k,
\]
then we have
\begin{equation}
  \label{eq:am-intro}
a_m = \frac{y_0^m(R_{k}^m g)(z_0)}{m!},  
\end{equation}
where $y_0 = \Im(z_0)$. Note that $R_{k}^m g$ is a non-holomorphic modular form of weight $k+2m$.

When $z_0 \in \Hb$ is a CM point, i.e.\ the unique root in $\Hb$ of a quadratic equation 
with integral coefficients
and discriminant $-D < 0$,  
the theory of complex multiplication implies that that the Laurent coefficients $a_m$ have good arithmetic properties.
If $g$ has algebraic Fourier coefficients and $z_0$ is a CM point of discriminant $-D < 0$, then the Laurent coefficient $a_m$ is an algebraic multiple of $\Omega_{-D}^{2m+k}$, 
where $\Omega_{-D}$ is the Chowla-Selberg period defined by
\begin{equation}
  \label{eq:Chowla-Selberg}
  \Omega_{-D} := \frac1{\sqrt{2\pi D}} \lp \prod_{j = 1}^{D-1} \Gamma \lp \frac{j}{D}\rp^{\chi_{-D}(j)} \rp^{w_{-D}/(4h_{-D})} .
\end{equation}
Here $\chi_{-D} (j) := (\frac{-D}{j})$ is the Dirichlet character, $w_{-D}$ and $h_{-D}$ are the number of roots of unity and class number of $\Qb(\sqrt{-D})$. 
For more details about Laurent expansions of modular forms we refer the reader to the exposition in Zagier's part of the book \cite{brdgza08}.

Unlike the Fourier coefficients, Rodriguez--Villegas and Zagier \cite{villegaszagier} showed that the Laurent coefficients $a_m$ can be obtained as special values of a family of polynomials defined by a ``quasi-recursion''. This allows for the quick numerical computation of the coefficients $a_m$, and is a useful tool for the study of Laurent expansions of holomorphic modular forms. For example, O'Sullivan and Risager \cite{osullivanrisager} used this method to compute the Laurent expansions of the discriminant $\Delta(z)$ at CM points, and they proved that these coefficients satisfy interesting periodicity properties modulo primes and are typically all non-vanishing. For instance, at the CM point $z_0 = i$ they found the Laurent expansion
\[
(1-w)^{-12}\Delta(z) = \sum_{m=0}^\infty\left(-4^3 \Delta(i)\left(\frac{-2\pi \Omega_{-4}^2
}{\sqrt{3}} \right)^m\right)\frac{p_m(0)}{m!}w^m, \qquad \left(w = \frac{z-i}{z+i}\right)
\]
where the polynomials $p_m(t) \in \Zb[t]$ are defined recursively by
\[
p_0 = 1, \quad p_1 = 0, \quad p_{m+1} = -2mtp_m + 6(t^2-1)p_m' - m(m+11)p_{m-1}.
\]
Since $i$ is an elliptic point of order $2$, the Laurent coefficients of $\Delta(z)$ with odd index $m$ vanish for trivial reasons. However, O'Sullivan and Risager showed that $p_{m}(0) \equiv 1 \pmod 5$ if $m \equiv 0 \pmod 4$ and $p_{m}(0) \equiv 3 \pmod 5$ if $m \equiv 2 \pmod 4$. In particular, the even index Laurent coefficients of $\Delta(z)$ at $z_0 = i$ are all non-vanishing. The analogous non-vanishing statement for the Fourier coefficients of $\Delta(z)$ is an open problem known as Lehmer's conjecture, and seems to be much deeper.
Similar periodicity properties modulo primes and non-vanishing results for the Laurent coefficients of the Jacobi theta function $\sum_{n \in \Zb}q^{n^2}$ were conjectured by Romik \cite{romik} and were proved shortly after by Scherer \cite{scherer} and Guerzhoy, Mertens and Rolen \cite{guerzhoymertensrolen}.

Another advantage of the Laurent expansion is that it can be used to study modular forms on co-compact groups, when the corresponding compact Shimura curve does not have cusps, so Fourier expansions are not available. Nelson \cite{nelson1,nelson2} and Voight and Willis \cite{voightwillis} have developed methods for the numerical computation of Laurent coefficients of modular forms on co-compact groups.
For a given co-compact group, it is possible to carry out the approach of Rodriguez--Villegas and Zagier (see \cite{BG12} for the Shimura curve $X^B(1)$ associated to the discriminant 6 quaternion algebra $B$).
However, there does not seem to be any systematic method to provably compute the Laurent coefficients in the  general co-compact case. 

The goal of this present paper is to give two different approaches of evaluating the Laurent coefficients
  of \emph{meromorphic} modular forms.
In the first approach, we extend the Rodriguez--Villegas--Zagier method in two ways to compute the Laurent coefficients of meromorphic modular forms on $\SL_2(\Zb)$ in terms of special values of recursively defined polynomials or rational functions. Method one, described in Proposition~\ref{prop:rec1}, involves a non-linear recursion for rational functions which works well in the case of a pole of small order. 
Method two, given in Theorem~\ref{thm:rec2}, involves a linear recursion of polynomials. To make the recursion start, one needs the Laurent expansion at the given CM point of the logarithmic derivative of a suitable modular function (this in turn can be easily obtained with the first method, since it has a pole of order one). The Laurent expansion of the log-derivative of the modular function has to be computed once and for all, and can be used to compute the Laurent expansion of any mermorphic modular form at the given CM point. 
\begin{theorem}
\label{thm:rec2}
Let~$g$ be a modular form of weight~$k$ on~$\SL_2(\Zb)$ with a pole of order~$N\ge0$ at a CM point~$z_0=x_0+y_0$ of discriminant $-D < 0$, and let~$h$ be a modular function with a simple zero in~$z_0$.
Denote $\sum_{m\ge-1}{h_m\Omega_{-D}^{2+2m}(-4\pi y_0w)^m}$  the Laurent expansion of~$h'/h$ at~$z_0$.

Then the Laurent expansion of~$g$ at~$z_0$ is given by
\[
(1-w)^{-k}g(z)=\sum_{m=-N}^\infty{q_m(t_0)\Omega_{-D}^{k+2m}(-4\pi y_0w)^m}
\] 
where~$t_0\in\bar{\Qb}$ and the polynomials~$q_m(t)$ satisfy the quasi-recursion
\begin{align*}
(m+1+N)q_{m+1}(t)&=a_1(t)q_m'(t) +a_2(m,k,t)q_m(t)+a_3(m-1,k,t)q_{m-1}(t)\\
&-N\sum_{j=-N}^{m}{(-1)^{m-j}q_j(t)\bigl(h_{m-j}+a_2(1,0,t)h_{m-j-1}+a_3(1,0,t)h_{m-j-2}}\bigr)
\end{align*}
for some polynomials~$a_i(x,y,t)$ with algebraic coefficients and initial data~$q_{-N}(t)$ dependent on~$g$ and~$h$.
\end{theorem}
\begin{remark}
For~$N=0$ this result reduces to the quasi-recursion of Rodriguez--Villegas and Zagier. 
\end{remark}

\begin{example}
If~$z_0=i$, the quasi-recursion given by Theorem~\ref{thm:rec2} is
\[
\begin{aligned}
(m+1+N)\cdot q_{m+1}(t)&=\frac{(t^2-1)}{2}q_m'(t)-\frac{(k+2m)t}{12}q_m(t)-\frac{(m+k-1)}{144}q_{m-1}(t)\\
&-N\sum_{j=-N}^m{(-1)^{m-j}q_j(t)\biggl(h_{m-j}+\frac{t}{6}h_{m-j-1}+\frac{h_{m-j-2}}{144}\biggr)}\,.
\end{aligned}
\]
One can choose~$h=E_6E_4^{-3/2}$. For the case~$N=0$ compare with~\eqref{eq:qrec}.
\end{example}

In the second approach, we give a new formula for the Laurent coefficients at CM points of meromorphic modular forms which arise as regularized theta lifts of weakly holomorphic modular forms of half-integral weight. For simplicity, we restrict to meromorphic modular forms for $\SL_2(\Zb)$ at certain CM points in the introduction. However, we would like to emphasize that the method works in much greater generality, including modular forms for congruence subgroups and co-compact groups (see the example concerning the Shimura curve $X^B(1)$ in Section \ref{subsec:ex2}).

Let $k \in 2\Zb$ be an even integer, and let $f = \sum_{n\gg -\infty}c_f(n)q^n \in M_{k+1/2}^{!,+}$ be a weakly holomorphic modular form of weight $k+1/2$ for $\Gamma_0(4)$ in the Kohnen plus space. For simplicity, we assume that $c_f(0) = 0$. Its Borcherds-Shimura lift is defined by a regularized theta integral and has the following Fourier expansion near the cusp infinity
\[
  \Phi_{\BS}(f,z) =
 - (-2)^{1 + k/2}
  \sum_{n\geq 1}\bigg(\sum_{d \mid n}\left(\frac{n}{d}\right)^{k-1}c_f(d^2) \bigg)q^n.
\]
It is a meromorphic modular form of weight $2k$ with poles at the CM points of discriminants $-D < 0$ for which $c_f(-D) \neq 0$.

Let $z_0 = \frac{-1 + \sqrt{-D}}2 \in \Hb$ be a CM point with  $-D < 0$ an odd, fundamental discriminant. 
Denote $\Oc := \Zb[z_0]$ the ring of integers of the imaginary quadratic field $K := \Qb(z_0)$. 
For any $k \in \Nb$ and $\mu \in \Zb/D\Zb$, we can associate a binary theta series
\begin{equation}
  \label{eq:thetaN}
  \theta_\mu^{(k)}(\tau) :=
  \sum_{\lambda \in \sqrt{-D}\Oc + \mu} \lambda^k q^{\Nm(\lambda)/D},
\end{equation}
which is a holomorphic cusp form of weight $k+1$.
Furthermore, it is the $\mu$-th component of a vector-valued modular form associated to the positive definite lattice $N = (\Oc, \Nm)$.
Now, let $\tilde\theta^{(k)}_\mu(\tau)$ be $\mu$-th component of a \textit{harmonic Maass form} with $\xi$-image $\theta_\mu^{(k)}$.
Denote $c_\mu^{(k),+}(\frac{n}D)$ the $(\frac{n}D)$-th Fourier of its holomorphic part, with $n$ necessarily integral.
In \cite{ehlenlischwagenscheidt}, we showed that these coefficients are in $K^{\mathrm{ab}} \cdot \Omega_{-D}^{2k}$. 
We now have the following result.

 \begin{theorem}
   \label{thm:intro}
      Let $f, z_0$ be as above, and 
	\[
	(1-w)^{-2k}\Phi_{\BS}(f,z) = \sum_{m  \gg -\infty}{a_m}{}w^m
	\]
	with $w = \frac{z-z_0}{z-\overline{z}_0}$ be the Laurent expansion of the Borcherds-Shimura lift $\Phi_{\BS}(f,z)$ around $z = z_0$. 
        If we  write $m = 2r + \ell \geq 0$ with $\ell \in \{0,1\}$, then we have
        \begin{equation}
          \label{eq:am-intro}
          \begin{split}
                      a_m =
\kappa
& \sum_{n, b \in \Zb} c^{(k+m), +}_{-b/2}\lp \frac{n}{D}\rp
  c_f\lp \frac{ -b^2 - 4n}D \rp b^\ell ( b^2 + 4n)^r P^{(\ell - 1/2, -k-m)}_r \lp \frac{4n - b^2 }{4n + b^2} \rp,
          \end{split}
        \end{equation}
where $\kappa = -\frac{  2^{3k/2 + 2m}}{{D^{k+m} \cdot m! \binom{-k -m + r}{r}}} \pi^m$ and
 $P^{(\alpha, \beta)}_r$ is the $r$-th Jacobi polynomial (see \eqref{eq:RCq}).
\end{theorem}

We refer to Theorem~\ref{thm:main} for the general statement and Section \ref{sec:4} for the proof of the above result.
\begin{example}
  \label{exmp:intro}
  Although the formula for $a_m$ looks rather involved on a first glimpse, it is easy to evaluate it numerically.
  For example, let $k = 2$ and
\begin{equation}
  \label{eq:fsc}
  f(\tau) =
  q^{-3} + 64q - 32384q^{4} + 131565q^{5} - 4257024q^{8} + 11535936q^{9} + O(q^{10}) \in M_{5/2}^{!, +}.
\end{equation}
Then around $z = z_0 = \frac{-1 + \sqrt{3}i}{2} \in \Hb$, we have
\begin{equation}
  \label{eq:ex1}
  \Phi_\BS(f, z) = - 2^8 \frac{\Delta(z)}{E_4^2(z)}
=  (1-w)^{4}
  \lp\frac{1}{9}\pi^{-2} w^{-2} - 2 \pi \Omega_{-3}^6 w + \dots \rp.
\end{equation}
For $m = 1$, we have $r = 0, \ell = 1$ and
$$
\theta^{(3)}_{\mu} = \mu \cdot  3\eta(\tau)^8
$$
for $\mu = 0, \pm 1$, where $\eta(\tau)$ is the Dedekind eta function.
In \cite[section 4]{ehlenlischwagenscheidt}, we constructed a $\xi$-preimage $\tilde\theta^{(3)}_\mu$ of $\theta^{(3)}_\mu$, whose holomorphic part  $\theta^{(3), +}_\mu$ has the Fourier expansion
$$
\theta^{(3), +}_{\mu} =  \frac{3^2 \Omega_{-3}^6}{2^{5}} \mu \cdot (q^{-1/3} - 61 q^{2/3} + ...).
$$
Theorem \ref{thm:intro} then gives the result
$$
a_1 =
-\frac{2^5 \cdot  \pi}{3^3} 
 \sum_{n, b \in \Zb} c^{(3), +}_{-b/2}\lp \frac{n}{3}\rp
  c_f\lp \frac{ -b^2 - 4n}3 \rp b = -\frac{\pi \Omega_{-3}^6}{3} \cdot (64 + 64 - 61 -61) = -2\pi \Omega_{-3}^6
  $$
  as only the terms $(n, b) = (-1, \pm 1), (2, \pm 1)$ contribute. 
\end{example}

\begin{remark}
  \begin{enumerate}
  \item 
    The main difficulty to apply Theorem \ref{thm:intro} is in finding an explicit $\xi$-preimage $\tilde\theta^{(k+m)}_{\mu}(\tau)$. In \cite{ehlenlischwagenscheidt}, we gave a general construction of such $\xi$-preimages for binary theta functions, also using Borcherds-Shimura lifts.
This has the advantage that the Fourier coefficients are provably algebraic with explicitly bounded denominators.
    However, in practice it is often simpler to write $\tilde\theta^{(k+m)}_{\mu}(\tau)$ as a linear combination of Maass-Poincar\'e series, whose coefficients can then be computed in terms of Kloosterman sums and Bessel functions. Combining with the main result in \cite{ehlenlischwagenscheidt}, one can then use these approximations to find the exact coefficients.
    \item 
In addition to the Borcherds-Shimura lift, there is a regularized theta lift from weakly holomorphic moludar forms of weight $3/2 - k$ to meromorphic modular forms of weight $2k$.
This was studied in detail by Zemel in \cite{Zemel15} and we call it the Borcherds-Zemel lift. As in the case of Borcherds-Shimura lift, its Laurent coefficients at CM points are related to Fourier coefficients of harmonic Maass
forms. The precise result is contained in Theorem \ref{thm:main2}.

\item
  We will give an example for the numerical evaluation of our formulas in Section~\ref{sec:4} in the case of Laurent expansion of a modular form on a Shimura curve.
\end{enumerate}
\end{remark}

The starting point of both approaches is the generalization of the expression \eqref{eq:am-intro} to the case when $z_0$ is a pole of $g$.
This is accomplished in Proposition \ref{prop:LaurentExp}, where the value of $R^mg$ at $z_0$ is replaced by the constant term of its Laurent expansion at $z = z_0$.
To prove Theorem \ref{thm:rec2}, one considers the generating function~$P(t,X)$ of the polynomials~$\{p_m(t)\}_{m\ge0}$ given by the Rodriguez--Villegas--Zagier quasi-recursion applied to the holomorphic modular form~$gh^N$. The ratio~$Q(t,X)$ of~$P(t,X)$ and (essentially) the~$N$-th power of the Laurent expansion of~$h$ is, after specialization of~$t$ to some~$t_0\in\bar{\Qb}$, the Laurent expansion of~$g$. A recursion for the polynomial coefficients of~$Q(t,X)$ is deduced by using a linear differential operator~$L_k$ that annihilates the generating series~$P(t,X)$.

To prove Theorem \ref{thm:intro} (or the more general version in Theorems \ref{thm:main} and \ref{thm:main2}), we express $R^m\Phi_\BS(f, z)$ as another regularized theta lift of $f$.
This is a consequence of the actions of differential operators on theta kernels.
Then at a CM point $z_0$, the signature $(1, 2)$ lattice $L$ contains the direct sum of a positive definite unary lattice $P$ and negative definite binary lattice $N^-$.
The theta kernel $\Theta_L$ is then essentially $\Theta_P \otimes\overline{\Theta_{N}}$. 
Using differential operators, such as the Rankin-Cohen bracket, we can construct
a preimage of  $\Theta_P \otimes\overline{\Theta_{N}}$ under the lowering operator out of 
$\Theta_P$ and the harmonic Maass form $\tilde\Theta_{N}$.
Finally an application of Stokes' Theorem expresses the value of $R^m\Phi_\BS(f, z)$ at $z = z_0$ as a linear combination of Fourier coefficients of $\tilde\Theta_{N}$.
This strategy has been employed often to study CM-values of regularized theta lifts. But it is the first time to be applied to study Laurent coefficients of modular forms at CM points.


The structure of the paper is as follows. In section \ref{sec:prelim}, we give preliminary information about modular forms, regularized theta lifts, differential operators, and Laurent expansions.
Most of the results are known, except possibly for Propositions \ref{prop:LaurentExp} and \ref{prop:over-reg}, which are needed for our purpose.
In section \ref{sec:rec}, we recall the recursion of Rodriguez--Villegas and Zagier, and generalize it to the setting of meromorphic modular forms of any weight.
In section \ref{sec:4}, we state and prove Theorems \ref{thm:main} and \ref{thm:main2} regarding the Laurent coefficients of regularized theta lifts.
The last subsection contains the example of the Laurent coefficients of a modular form on a compact Shimura curve. 

\subsection*{Acknowledments}
We thank Brandon Williams for his help with using and modifying his wonderful code \cite{WeilRep} in computing vector-valued modular forms. 
G.\ Bogo and Y. Li are supported by the Deutsche Forschungsgemeinschaft (DFG) through the Collaborative Research Centre TRR 326 ``Geometry and Arithmetic of Uniformized Structures", project number 444845124. 
M. Schwagenscheidt was supported by SNF project PZ00P2\_202210.

\section{Preliminaries}
\label{sec:prelim}
\subsection{Modular forms for the Weil representation}

Let $V$ be a rational quadratic space of signature $(b^+,b^-)$ with quadratic form $Q : V \to \Qb$ and corresponding bilinear form $(\cdot,\cdot)$.
Denote $G$ the metaplectic cover of $\SL_2$, $H_V = \GSpin(V)$ and $\omega = \omega_\psi$ the Weil representation of $G(\Ab) \times H_V(\Ab)$ on the space of Schwartz functions $\Sc(V(\Ab))$ with $\psi$  the standard additive character on $\Ab/\Qb$.

Let $L \subset V$ be an even lattice and let $L'$ be the dual lattice of $L$. We also let $L^-$ (resp.\ $V^-$) denote the lattice $L$ (resp.\ vector space $V$) with quadratic form $-Q$. The discriminant group $L'/L$ is finite and isomorphic to $\hat{L}'/\hat{L}$, where $\hat{L}' := L' \otimes \hat{\Zb}$. For $\mu \in L'/L$ we consider the characteristic functions 
\begin{equation}
  \label{eq:phimu}
  \phi_\mu := \cha(\hat L + \mu), \quad \mu \in L'/L,
\end{equation}
They form a basis for the subspace
\begin{equation}
  \label{eq:SL}
  \Sc_L := \bigoplus_{\mu \in L'/L}\Cb \phi_\mu \subset \Sc(\hat V)
=  \bigcup_{L \subset V \text{ lattice}} \Sc_L  
\end{equation}
of Schwartz functions on $\hat V := V(\hat\Qb)$ that are supported on $\hat L'$ and constant on cosets of $\hat{L} := L \otimes \hat{\Zb}$.
For a sublattice $L_0 \subset L$, we have $L' \subset L_0'$ and a natural 
inclusion $\Sc_{L} \subset \Sc_{L_0}$.
We let $\langle \cdot,\cdot \rangle$ be the bilinear pairing on $\Sc_L$ defined by $\langle \phi_\mu,\phi_\nu \rangle = \delta_{\mu,\nu}$.
More generally for $\phi_1, \phi_2 \in \Sc(\hat V)$ with $\phi_i \in \Sc_{L_i}$, we view $\phi_i \in \Sc_L$ with $L = L_1 \cap L_2$ and can define $\langle \phi_1, \phi_2\rangle$ accordingly.

Let $\rho_L$ be the Weil representation of the two-fold metaplectic cover $\Mp_2(\Zb)$ of $\SL_2(\Zb)$ on $\Sc_L$ (see \cite[Section~4]{borcherds}).
For $k \in \frac{1}{2}\Zb$, let $H_{k, L}$ be the space of harmonic (weak) Maass forms of weight $k$ valued in $\Sc_L$ for the representation $\rho_L$  (see \cite{bruinierfunke04}) for which there is a Laurent polynomial 
\[
P_f(\tau) = \sum_{\mu \in L'/L}\sum_{\substack{m \in \Qb_{\leq 0} \\ m \gg -\infty}}c_f^+(m,\mu)q^m \phi_\mu,
\]
called the principal part of $f$, such that $f(\tau)-P_f(\tau) = O(e^{-\varepsilon})$ as $v=\Im(\tau) \to \infty$, for some $\varepsilon > 0$. It contains $S_{k,L} \subset M_{k,L} \subset M_{k,L}^!$,  the subspaces of cusp forms, holomorphic, and weakly holomorphic modular forms.
For a subgroup $N$ (resp.\ a character $\chi$) of the finite orthogonal group $\SO(L'/L)$, we denote $H^{N}_{k, L}$ (resp.\ $H^\chi_{k, L}$) to denote the $N$-invariant (resp.\ $\chi$-isotypic) subspace of $H_{k, L}$.

When 
$
\kappa := k-\frac{b^--b^+}{2}
$
is an integer, the components of $f(\tau) = \sum_{\mu \in L'/L}f_\mu(\tau)\phi_\mu \in H_{k, L}$ satisfy the symmetry
\begin{equation}
  \label{eq:sym}
f_{-\mu}(\tau) = (-1)^\kappa f_{\mu}(\tau)  .
\end{equation}
Otherwise, the space $H_{k,L}$ is trivial.
Let $f^+$ denote the holomorphic part of $f$ with Fourier expansions of the shape
\begin{align*}
f^+(\tau) &= \sum_{\mu \in L'/L}\sum_{\substack{m \in \Qb \\ m \gg -\infty}}c_f^+(m,\mu)q^m \phi_\mu,~ c_{f}^\pm(m,\mu) \in \Cb.
\end{align*}
The $\xi$-operator
\[
\xi_{k}f = 2iv^k \overline{\partial_{ \overline{\tau}}f}
\]
defines a surjective map $\xi_k : H_{k,L} \to S_{2-k,L^-}$.

For $m \in \Qb$, denote
\begin{equation}
  \label{eq:varphim}
  \varphi_{f,m} := \sum_{\mu \in L'/L} c_f^+(m, \mu) \phi_\mu \in \Sc_L. 
\end{equation}
which is non-trivial only for finitely many negative $m$.
By \eqref{eq:sym}, we have
\begin{equation}
  \label{eq:parity}
  \varphi_{f,m}(-x) = (-1)^{\kappa}  \varphi_{f,m}(x)
\end{equation}
for all $x \in \hat V$.

\subsection{Differential operators}
For $k \in \frac{1}{2}\Zb$, we have
the raising and lowering operator
\begin{equation} \label{defops}
  \begin{split} R_{\tau, k}&= R_{k} := 2i\partial_{\tau}+\frac{k}{v}, \quad L=L_{\tau}:=-2iv^{2}\partial_{\overline{\tau}},
  \end{split}
\end{equation}
Note that $\xi_k f = v^{k-2} \overline{L_\tau f}$. 
For $n \in \Nb$, we write
\begin{equation}
  \label{eq:Rn}
R^n =  R^n_k := R_{k+2n-2} \circ R_{k+2n-4} \circ \dots \circ R_{k+2} \circ R_{k}  
\end{equation}
for the iterated raising operator.

For modular forms $f_i$ of weight $k_i \in\frac{1}{2}\Zb$ and $r \in \Nb$, we have the usual Rankin-Cohen operator
\begin{equation}
  \label{eq:RC}
  \begin{split}
    [f_1, f_2]_r(\tau) &:= (2\pi i)^{-r} \sum^r_{s = 0} (-1)^s \binom{k_1 + r - 1}{s} \binom{k_2 + r - 1}{r - s}    (\partial_{\tau}^{(r-s)}f_1) (\tau)(\partial_{\tau}^{(s)} f_2) (\tau) \\
    &= (4\pi)^{-r} \sum^r_{s = 0} (-1)^{r-s} \binom{k_1 + r - 1}{s} \binom{k_2 + r - 1}{r - s}
    (R^{r-s}_{} f_1)(\tau)(  R^s_{} f_2)(\tau)
  \end{split}
\end{equation}
where
$\binom{m}{n} := \frac{m(m-1)(m-2)\dots(m-n + 1)}{n!}$ is the binomial coefficient.
The result $[f_1, f_2]_r$ is modular of weight $k_1 + k_2 + 2r$. 
Furthermore, suppose $f_1$ is harmonic and $f_2$ is holomorphic, then
\begin{equation}
  \label{eq:RCid}
  \begin{split}
 (4\pi)^{r} \binom{k_1 +r -1}{r}^{-1}      L_\tau [f_1, f_2]_r
 & =
    (L_\tau f_1)( R^r_{} f_2)
   =    R^r_{}L_\tau (f_1f_2)
  \end{split}
\end{equation}
for any $r \in \Nb$. 
More generally for two vector-valued real-analytic functions $f = \sum_\mu f_\mu \vf_\mu$ and $g = \sum_{\nu} g_\nu \wf_\nu$, their Rankin-Cohen bracket is defined componentwisely as
\begin{equation}
  \label{eq:vvRC}
  [f, g]_r = [f, g]^{(k_1, k_2)}_r := \sum_{\mu, \nu} [f_\mu, g_\nu]_r \vf_\mu \otimes \wf_\nu. 
\end{equation}
In particular for $f = q^n, g = q^m$ with $m + n \neq 0$, we have
\begin{equation}
  \label{eq:RCq}
  [f, g]_r = (m+n)^r P_r^{(k_2 - 1, k_1 - 1)}\lp \frac{n-m}{n+m} \rp q^{m+n},
\end{equation}
where 
\[
P_r^{(\alpha, \beta)}(x) = \sum_{s = 0}^r \binom{r+\alpha}{r-s} \binom{r+\beta}{s} \lp \frac{x+1}2 \rp^{r-s} \lp \frac{x - 1}2 \rp^s
\]
is the Jacobi polynomial. Using the expression 
\[
P_r^{(\alpha, \beta)}(x)
= \frac{(\alpha + 1)_n}{n!} \prescript{}{2}F_1\left(-r, 1 + \alpha + \beta + r; \alpha + 1; \frac12(1-z)\right)
\]
and the Pfaff transformation 
\[
\prescript{}{2}F_1(a, b; c; z) = (1-z)^{-a} \prescript{}{2}F_1\left(a, c-b; c; \frac{z}{z-1}\right),
\]
one can show that
\begin{equation}
  \label{eq:CT-RC}
  \CT([f, [g, h]_r]_0) =   \CT([[f, g]_r, h]_0)  =   \CT([[h, f]_r, g]_0)
\end{equation}
for Fourier series $f, g, h$ with weights $k_1, k_2, k_3$ satisfying $k_1 + k_2 + k_3 + 2r = 2$.

\subsection{Laurent expansions of meromorphic modular forms}
At $z_0=x_0+iy_0\in\mathcal{H}$, we take the local coordinate 
\begin{equation} \label{Awz}
  w
  :=\frac{z-z_0}{z-\overline{z_0}}\in\big\{w\in\mathbb{C} \, : \, |w|<1\big\}.
\end{equation}
It is easy to check that
\begin{equation}
  \label{eq:yid}
  z = \frac{w\bar z_0 - z_0}{w - 1}, \quad
  y =  \frac{1 - w \bar w}{(1-w)(1-\bar w)}y_0.
\end{equation}

For  a real-analytic function $f$,  we denote
\begin{equation}
  \label{eq:tf}
  \tf(w) := f \lp \frac{z_0 - w \overline{z}_0}{1-w} \rp ,
\end{equation}
which  has the following Laurent expansion at $z_0$,
  $$
  \tf(w) 
=  \sum_{m, n \gg -\infty} a^{}_{m, n} w^m (\bar w)^n
$$
with $a^{}_{m, n} \in \Cb$.
Clearly $f$ is meromorphic at $z_0$ if and only if $a^{}_{m, n} = 0$ for all $n \neq 0$, in which case we write $a_m := a_{m, 0}$. 
We say that $f$ has a \textit{meromorphic singularity at $z_0$} if $a^{}_{m, n} = 0$ for all $n < 0$.
This condition is preserved by the raising operator. 

When $f$ is a real-analytic modular form of weight $\kappa \in \Zb$, it is more natural to consider the Laurent expansion of $(1-w)^{-\kappa}\tf(w)$ \cite{brdgza08} and we denote
\begin{equation}
  \label{eq:CTk}
  \CT_{z_0}^{(\kappa)}(f) :=
  \CT\left((1-w)^{-\kappa}\tf\right).
\end{equation}
When $z_0$ is fixed, we will omit it from the notation.
For a modular form $f$, the weight is also inherent in it and we omit $\kappa$ from the notation as well.
\begin{lemma}
  \label{lemma:CT}
  Let $f$ be a real-analytic function.
  Then for any $\kappa \in \Zb$, we have
  \begin{equation}
    \label{eq:CTR}
    \partial_w ((1-w)^{-\kappa} \tf(w))
         = y_0 (1-w)^{-\kappa-2} \widetilde{R_\kappa f} (w)
      + \frac{\kappa \tf(w)}{(1-w)^{\kappa + 1}}
\frac{      1 + w}{1 + w \bar w} \bar w.
  \end{equation}
\end{lemma}

\begin{proof}
  This follows from the straightforward calculation
  \begin{align*}
    \partial_w ((1-w)^{-\kappa} \tf(w))
    &= (1-w)^{-\kappa-2} \lp (1- w) \kappa \tf(w)
      +  2iy_0 \partial_z f (z) \rp\\
    &     = y_0 (1-w)^{-\kappa-2} \widetilde{R_\kappa f} (w)
      + \frac{\kappa \tf(w)}{(1-w)^{\kappa + 2}}
      ((1 - w) - y_0/y)
  \end{align*}
  together with \eqref{eq:yid}.
\end{proof}
We now have the following result, which slightly generalizes Proposition 17 in Zagier's part of \cite{brdgza08}.
\begin{proposition}
  \label{prop:LaurentExp}
  Let $f$ be a mermorphic modular form of weight $\kappa \in \Zb$ with the Laurent expansion
  $$
  (1-w)^{-\kappa} \tf(w) = \sum_{m \gg -\infty} a_m w^m.
$$
Then for $m \ge 0$, we have
\begin{equation}
  \label{eq:LaurentC}
a_m = \frac{y_0^m}{m!} \CT(R_\kappa^m f)
\end{equation}
\end{proposition}

\begin{proof}
From \eqref{eq:CTR}, we have
  \begin{align*}
    (\partial_w^m (( 1- w)^{-\kappa} \tf(w)))
    &= \partial_w^{m-1}
    \lp    y_0 (1-w)^{-\kappa-2}
\widetilde{R_\kappa f} (w)
      +
\frac{\kappa \tf(w)}{(1-w)^{\kappa + 1}}
\frac{      1 + w}{1 + w \bar w} \bar w
      \rp\\
    &=         y_0^m (1-w)^{-\kappa-2m}
\widetilde{R^m_\kappa f} (w)
      +
\sum_{j = 1}^m
      \partial_w^{m-j}
      \lp \frac{\kappa \widetilde{R_\kappa^{j-1} f}(w)}{(1-w)^{\kappa + 1}}
      \frac{      1 + w}{1 + w \bar w} \rp
      \bar w.
  \end{align*}
  Since $f$ is meromorphic, the function $R^j_\kappa f$ has a meromorphic singularity at $z_0$ for all $j \ge 0$.
  Taking constant terms of both sides proves the result.
\end{proof}
For holomorphic modular forms an effective way to compute the Laurent coefficients was given in~\cite{villegaszagier}. We will recall this in Section~\ref{sec:rec}, and consider the case of Laurent expansions at poles. 

\subsection{Quadratic space and CM points.}
\label{subsec:qs}
Let $(V, Q)$ be a rational quadratic space of signature $(1, 2)$ and $H = H_V = \GSpin(V)$.
Throughout we fix an isometry
\begin{equation}
  \label{eq:isometry}
V(\Rb) \cong M_2(\Rb)^0 = \left\{\begin{pmatrix}b/2 & c \\ -a & -b/2\end{pmatrix} \in M_2(\Rb): a, b, c \in \Rb\right\}  
\end{equation}
with determinant as the quadratic form.
Furthermore, identify 
\[
(M_2(\Rb)^0, \det) \cong (\Rb^{1, 2}, (x_1^2 - x_2^2 - x_3^2)/2)
\]
with the basis $\{w_1, w_2, w_3\}$ given by
\begin{equation}
  \label{eq:vs}
  w_1 := \frac{1}{\sqrt{2}} \pmat{0}{-1}{1}{0}, \quad
  w_2 :=  \frac{1}{\sqrt{2}} \pmat{0}{1}{1}{0}, \quad
  w_3 :=  \frac{1}{\sqrt{2}} \pmat{1}{0}{0}{-1}.
\end{equation}
This gives us $\Sc(\Rb^{1, 2}) \cong \Sc(V(\Rb))$ and 
\begin{equation}
  \label{eq:bases-change}
  \pmat{b/2}{c}{-a}{-b/2} = \sum_{i = 1}^3 x_i w_i, \quad
x_1 =  \frac{-a - c}{\sqrt{2}},~
x_2 =  \frac{-a+c}{\sqrt{2}},~
x_3 =  \frac{b}{\sqrt{2}}.
\end{equation}

The group $H(\Rb) \cong \mathrm{GL}_2(\Rb)$ acts on $V(\Rb)$ via the conjugation action.
Denote  $Z_0 :=  \Rb X_0 \oplus \Rb Y_0$ with
\begin{equation}
  \label{eq:basis}
  X_0 := \begin{pmatrix}0 & 1 \\ 1 & 0\end{pmatrix}   ,\quad
  Y_0 := \begin{pmatrix}1 & 0 \\ 0 & -1 \end{pmatrix}, \quad
Z_0^\perp := \begin{pmatrix}0 & -1 \\ 1 & 0 \end{pmatrix},
\end{equation}
an ordered orthogonal basis of $M_2(\Rb)^0$.
When $V$ is isotropic, we fix a $\Qb$-isometry $(V, Q) \cong (M_2(\Qb)^0, N \cdot \det)$ with $N \in \Nb$. Composing this with  $(M_2(\Qb)^0, N \cdot \det) \cong (\sqrt{N} M_2(\Qb)^0, \det)$ and tensoring with $\Rb$ gives us an isometry $V(\Rb) \cong M_2(\Rb)^0$, which we require to be the same as in \eqref{eq:isometry}.

  For $z = x +iy \in \Hb$, define
  \[
h_z := n(x)m(\sqrt{y}) \in \SL_2(\Rb) \subset \GL_2(\Rb) \cong H(\Rb)
\]
with $n(b) := \smat{1}{b}{}{1}, m(a) := \smat{a}{}{}{a^{-1}}$, which  satisfies  $h_z \cdot i = z$. 
Let $\Db$ be the Grassmannian of oriented negative two-planes in $V(\Rb)$.
It has two connected components $\Db^\pm$ with $\Db^+$ containing $Z_0$. 
We can identify $\Hb$ with $\Db^+$ via
\begin{align}\label{identification DH}
  \begin{split}
      z &\mapsto h_z \cdot Z_0
      =       \Rb \Re Z(z) \oplus \Rb \Im Z(z) = (Z^\perp(z))^\perp,\\
      Z(z) &:= \begin{pmatrix}z & -z^2 \\ 1 & -z \end{pmatrix} = y(h_z \cdot (X_0 + iY_0)), \quad  Z^\perp(z) := y^{-1} \begin{pmatrix}x & -|z|^2 \\ 1 & -x \end{pmatrix} = h_z \cdot Z^\perp_0.
  \end{split}
\end{align}
Note that for any $\gamma = \smat{*}{*}{c}{d} \in \SL_2(\Rb)$, we have 
\begin{equation}
  \label{eq:decompose}
  h_{\gamma \cdot z} = \gamma h_z \kappa(\theta(\gamma, z)), \quad \kappa(\theta) := \smat{\cos\theta}{-\sin\theta}{\cos\theta}{\sin\theta}, \quad
  e^{i \theta(\gamma, z)} = |cz + d|/(cz+d).
\end{equation}
For convenience, we record the following change of variables
\begin{equation}
  \label{eq:change-var}
  \begin{split}
    (Z(z), Z(z_0))
    &= - 4y_0^2 \frac{w^2}{(1-w)^2},\quad
      (Z(z), Z(\overline{z_0}))
      = - 4y_0^2 \frac{1}{(1-w)^2},\quad\\
    (Z(z), Z^\perp(z_0))
    &= - 4y_0 \frac{w}{(1-w)^2},\quad
    (Z^\perp(z), Z({z_0}))
    = 4y_0 \frac{ w}{1-|w|^2},\\
    (Z^\perp(z), Z(\overline{z_0}))
    &= 4y_0 \frac{\bar w}{1-|w|^2},\quad\quad
      (Z^\perp(z), Z^\perp({z_0}))
    = 2 \frac{1 + |w|^2}{1 - |w|^2},\\
  \end{split}
\end{equation}
where $z = x + i y, z_0 = x_0 + iy_0$ and $w = \frac{z-z_0}{z-\overline{z}_0}$ as in \eqref{Awz}. 
For any compact open subgroup $K \subset H(\Ab_f)$, the associated Shimura variety
\begin{equation}
  \label{eq:XK}
X_K := H(\Qb) \backslash \Db \times H(\Ab_f)/K
\end{equation}
is a disjoint union of modular (resp.\ Shimura) curves when $V$ is isotropic (resp.\ anisotropic).  

We now quickly recall CM points and special divisors on $X_K$ (see section 2 of \cite{BKY12}). 
For a negative definite subspace $U^- \subset V$, we have two points $\{z_U^\pm \} \subset \Db$ given by $U^-(\Rb) \subset V(\Rb)$ with two possible orientations.
The group $T = \GSpin(U^-) = \GSpin(U)$ is isomorphic to $\mathrm{Res}_{E/\Qb} \Gm$ for an imaginary quadratic field $E$, and embeds into $H$, which gives us a CM cycle
\begin{equation}
  \label{eq:ZU}
Z(U, h) = T(\Qb) \backslash (\{z_U^\pm \}\times T(\hat \Qb)/(hKh^{-1} \cap T(\hat \Qb))) \to X_K,\quad [z_U^\pm, t] \mapsto [z_U^\pm, th]
\end{equation}
for any $h \in H(\hat\Qb)$.

In the signature $(1, 2)$ case here, the CM points coincide with special divisors defined in the following way.
For $w \in V(\Qb)$ with positive norm, its stabilizer $H_w$ in $H$ is isomorphic to $T= \GSpin(U)$ with $U^- := w^\perp$.
The analytic divisor $\Db_w = \{z \in \Db: z \perp w\} \subset \Db$ is simply $\{z_U^\pm\}$, from which we can construct a divisor $Z(w, h)$ on $X_K$ for any $h \in H(\hat\Qb)$, which is the image of
$$
H_w(\Qb) \backslash \Db_w \times H_w(\hat\Qb)/(H_w(\hat\Qb) \cap h K h^{-1})
\to X_K,~ (z, h_1) \mapsto (z, h_1h).
$$
Clearly, we have $Z(w, h) = Z(U, h)$. 
For $\phi \in \Sc(\hat V)^K$ and $m \in \Qb_{> 0}$ such that there exists $w \in V(\Qb)$ with $Q(w) = m$, we can define the weighted cycle
$$
Z(m, \phi) := \sum_{h \in H_w(\hat\Qb)\backslash H(\hat\Qb)/K} \phi(h^{-1} w) Z(w, h).
$$
When $\phi = \phi_\mu$ as in \eqref{eq:phimu}, we denote $Z(m, \mu) := Z(m, \phi_\mu)$.

\subsection{Theta functions}
In this section, we recall unary and binary theta functions and show that the theta kernel realizing the Shimura correspondence splits into products of these theta functions at CM points.

Let $(W, Q_W)$ be a positive definite, 1-dimensional quadratic space.
Fix a vector $w_0 \in W(\Rb)$ such that $Q_W(w_0) = 1/2$. 
Then the map  $w \mapsto (w, w_0)$ gives an isometry $W(\Rb) \cong (\Rb, x^2/2)$.
For $\ell \in \Nb_0$ and $\varphi_f \in \Sc(\hat W)$, we have the definite theta function
\begin{equation}
  \label{eq:thetaW}
  \begin{split}
    \theta_W^{(\ell)}(\tau, \varphi_f)
    &:=
\sqrt{v}^{-\ell - \frac{1}{2}}
\sum_{w \in W(\Qb)}
(\omega_W(g_\tau) \varphi^\ell)(w)\\
&=
(2\sqrt{ v})^{-\ell}
\sum_{w \in W(\Qb)}\varphi_f(x)
 \He_{\ell}(\sqrt{2\pi v}(w, w_0))q^{Q_W(w)},
  \end{split}
\end{equation}
where $\varphi^\ell = \varphi_f \varphi_\infty^\ell$ with  $W(\Rb) \cong \Rb$,
$$
\varphi_\infty^\ell(x) =
(-2\sqrt{2\pi})^{-\ell}(\partial_x  - 2\pi x)^\ell e^{-\pi x^2}
=
2^{-\ell}\He_{\ell}(\sqrt{2\pi}x)e^{-\pi x^2}
\in \Sc(\Rb)
\cong \Sc(W(\Rb)),
$$
and 
\[
\He_\ell(\xi) := (-1)^{\ell}e^{\xi^{2}}\big(\tfrac{d}{d\xi}\big)^{\ell}e^{-\xi^{2}}= \big(2\xi-\tfrac{d}{d\xi}\big)^{\ell}\cdot1
\]
is the $\ell$-th (physicist's) Hermite polynomial, 
and $\omega_W$ is the Weil representation associated to $W$.
It is a real-analytic modular form of weight $\ell + \frac{1}{2}$. 
When $\ell = 0, 1$, it is furthermore holomorphic. 
For an even lattice $P \subset W$, we denote
\begin{equation}
  \label{eq:ThetaP}
  \Theta_P^{(\ell)}(\tau) := \sum_{\mu \in P'/P} \theta^{(\ell)}_W(\tau, \phi_\mu),
\end{equation}
with $\phi_\mu := \cha(\hat P + \mu) \in \Sc(\hat W)$. 


For a positive definite quadratic space $(U, Q_U)$ of dimension 2, there exists
$a \in \Qb_{>0}$ and an imaginary quadratic field $E$ such that
\begin{equation}
  \label{eq:Ua}
(U, Q) \cong (E, a \Nm)  .
\end{equation}
Tensoring this with $\Rb$ and composing with $E \otimes\Rb \cong \Cb, z \mapsto \sqrt{a} z$, we have the following isometry of real quadratic spaces. 
\begin{equation}
  \label{eq:UR}
(  U(\Rb), Q)\cong (\Cb, \Nm_{\Cb/\Rb}).
\end{equation}
We identify $T := \GSpin(U) \cong\mathrm{Res}_{E/\Qb} \Gm$, and
\begin{equation}
  \label{eq:SOU}
  \SO_U(R) \cong R^1 := \k^1 \otimes R, \quad \k^1 := \{z \in \k: \Nm(z) = 1\},
\end{equation}
for any $\Qb$-algebra $R$.
Under these isomorphisms, the surjective map $T \to \SO_U$ is given by
$$
\Nm^-: T \to \SO_U, \quad
z \mapsto  z/\bar{z},
$$
with $\bar{z}$ the Galois conjugation of $z$.
 Let $\phi^k = \phi_f \phi^k_\infty \in \Sc(U(\Ab))$ be a Schwartz function such that $\phi_f \in \Sc(\hat U)  \cong \Sc(\widehat{ U^-})$ and
$$
\phi^k_\infty: U( \Rb) \cong \Cb \to \Cb, \quad  z\mapsto  z^{k} e^{-2\pi z \overline z}.
$$
Then we define the binary theta series
\begin{equation}
  \label{eq:thetaUk}
      \theta_{U}^{(k)}(\tau, h, \phi_f)
      := 
      \sqrt{v}^{-k-1}\sum_{x \in U} (\omega_U(g_\tau, h) \phi^k)(x)
      = \sqrt{a}^k \sum_{x \in E} x^k \phi_f(h^{-1} x) \ebf(a \Nm(x))
\end{equation}
with $h \in T(\hat\Qb)$, which is a holomorphic modular form in $\tau$ of weight $2k+2$.
Note that this theta function depends on the isometry in \eqref{eq:Ua}. 
For an even, integral lattice $N \subset U$, we denote
\begin{equation}
  \label{eq:ThetaN}
  \Theta_N^{(k)}(\tau, h) := \sum_{\mu \in N'/N} \theta^{(k)}_U(\tau, h, \phi_\mu),
\end{equation}
with $\phi_\mu := \cha(\hat N + \mu) \in \Sc(\hat U)$.
We omit $h$ from the notation when it is 1.

Let $V$ be a quadratic space of signature $(1, 2)$ as in section \ref{subsec:qs}.
and define the theta function
  $$
  \theta_{V}(\tau, (z, h), \varphi) = 
  \sum_{x \in V(\Qb)}
(\omega_V(g_\tau, h_z) \varphi)(h^{-1}x)
$$
for $\varphi \in \Sc(V(\Ab))$ and $(z, h) \in \Db \times H(\hat\Qb)$.
When $\varphi = \varphi^{(\ell, k)} = \varphi_f \varphi_{\infty}^{(\ell, k)} \in \Sc(V(\Ab))$ with $\ell, k \in \Nb$
and
\begin{equation}
  \label{eq:varphik}
  \varphi_{\infty}^{(\ell, k)}(x_1, x_2, x_3) :=
    2^{-\ell}\He_\ell( \sqrt{2\pi} x_1)
    \lp x_2 - ix_3 \rp^{k}
    e^{-\pi (x_1^2  + x_2^2 + x_3^2 )}
  \end{equation}
  a Schwartz function in $ \Sc(\Rb^{1, 2})$, 
we denote
\begin{equation}
  \label{eq:thetak}
  \theta_{V}^{(\ell, k)}(\tau, (z, h), \varphi_f) := \sqrt{v}^{k  - \ell + \frac{1}{2}} y^{-k}\theta(g_\tau, (z, h), \varphi^{(\ell, k)})
\end{equation}
which is a real-analytic modular form of weight $2k$ in on $X_K$ with $K$ any open compact that preserves $\varphi_f$.
In the variable $\tau$, it is modular of weight $\ell - k - 1/2$.

For an even, integral lattice $L \subset V$, 
let $K(L) \subset H(\Ab_f)$ be the largest open compact subgroup preserving $\hat L$ and acting trivially on $\hat L'/\hat L \cong L'/L$.
If $K$ is contained in $K(L)$, then we denote
\begin{equation}
  \label{eq:thetaLk}
  \begin{split}
    \Theta_{L}^{(\ell, k)} (\tau, z, h) &:= \sum_{\mu \in L'/L} \theta_V^{(\ell, k)} (\tau, (z, h), \phi_\mu)  \phi_\mu\\
&=
    v^{k+1 - \ell/2} y^{-2k}
    \sum_{\substack{\mu \in L'/L\\\lambda \in h(L+\mu)}}
    2^{-\ell}\He_\ell \lp \sqrt{\pi v}
    \lp \lambda, Z^\perp(z)\rp\rp\\
&\times  \left(\lambda,
  -\frac{1}{\sqrt{2}}
  Z(\bar z)
  \right)^{k}
  \ebf(Q(\lambda_{z^\perp})\tau + Q(\lambda_{z})\overline{\tau})\phi_{\mu}
  \end{split}
\end{equation}
 the vector-valued theta kernel with $\phi_\mu$ as in \eqref{eq:phimu} with $[z, h] \in X_K$. Here $\lambda_z$ denotes the orthogonal projection of $\lambda$ to the negative plane corresponding $z$, and $\lambda_{z^\perp}$ denotes the projection of $\lambda$ to the orthogonal complement in $V(\Rb)$ of the plane corresponding to $z$.

 Let $U^- \subset V$ be a negative definite 2-dimensional subspace and denote $W := (U^-)^\perp$.
 Then $U^-$ gives rise to a CM cycle $Z(U)$, which consists of toric orbits of the CM point $z_U^\pm$ coming from $U^- \otimes \Rb$.
 Denote $z_U = x_U + iy_U := z_U^+$ the point in the upper half plane. 
 There are unique square-free $a, M \in \Nb$ such that
 \begin{equation}
   \label{eq:WUisom}
   \begin{split}
     W & \cong (\sqrt{M} \Qb, x^2),~ \lambda \mapsto  (\lambda, Z^\perp(z_U)),\\
     U &\cong (\sqrt{a} E, \Nm),~  \lambda \mapsto  -(2y_U)^{-1}(\lambda, Z({z_U})).
   \end{split}
 \end{equation}
 Then there is a canonical isomorphism
\begin{align*}
\iota:  \Sc(\hat{U}) \otimes \Sc(\hat W)  &\to \Sc(\hat V) ,\\
\varphi_U\otimes \varphi_W &\mapsto \varphi: x = (x_U, x_W) \mapsto \varphi_U(x_U) \varphi_W(x_W).
\end{align*}
At a CM point $[z_U^\pm, h_0] \in Z(U)$, the theta function $\theta_V^{(\ell, k)}$ splits as
\begin{equation}
  \label{eq:theta-split}
  \theta_{V}^{(\ell, k)}(\tau, [z_U^\pm, h_0], \varphi_f)
  =
(\sqrt{2}/y_U)^{k} \theta_W^{(\ell)}(\tau, \varphi_W)
\overline{  \theta_{U}^{(k)}(\tau, h_0, \phi_f) } v^{k+1}
\end{equation}
when $\varphi_f = \iota(\phi_f \otimes \varphi_W)$.

At a CM point $[z_U^\pm, h] \in Z(U) \subset X_K$, denote $N^- := L \cap U^-$ and $P := L \cap W$. Then $P \oplus N^-$ is a sublattice of $L$ and the difference
$$
\Theta_L^{(\ell, k)}(\tau, z_U^\pm, h) 
- (\sqrt{2}/y_U)^{k} \Theta^{(\ell)}_P(\tau) \otimes \overline{\Theta^{(k)}_N(\tau, h)} v^{k+1} 
$$
as a function in $\Sc(\hat V)$ is identically zero on $\hat L' \subset \hat V$.
Therefore for any $\phi \in \Sc(\hat V)$ with support in $\hat L'$, we have
\begin{equation}
  \label{eq:sub-pair}
\langle \phi, \Theta_L^{(\ell, k)}(\tau, z_U^\pm, h)  \rangle
= (\sqrt{2}/y_U)^{k} \langle \phi, \Theta^\ell_P(\tau) \otimes \overline{\Theta^k_N(\tau, h)} v^{k+1} \rangle.
\end{equation}
We remark that in the literature identities of this shape are often stated using a ``restriction operator'' $\phi \mapsto \phi_{P \oplus N}$ from modular forms for $\rho_L$ to modular forms for $\rho_{P \oplus N}$ (see for example the proof of Lemma~4.4 in \cite{bruinieryang}), that is, $\phi$ on the right-hand side is replaced with $\phi_{P \oplus N}$. However, due to our definition of the bilinar pairing $\langle \cdot, \cdot \rangle$, this is not necessary here.

\subsection{Regularized theta lifts}
Let $F(\tau)$ be a smooth function on $\Hb$ which is invariant with respect to $\tau \mapsto \tau + 1$. 
Suppose that the constant term of its Fourier expansion has at most polynomial growth as  $v = \Im(\tau)$ goes to infinity.
Let  $\Fc_T \subset \Hb$ be the truncated fundamental domain of $\SL_2(\Zb)\backslash \Hb$ at height $T$.
For $\Re(s) \gg 0$, the limit
$\lim_{T\to\infty}\int_{\Fc_T} F(\tau) v^{-s} d\mu(\tau)$ exists, and defines a holomorphic function in $s$ that can be analytically continued to $s \in \Cb$.
We can define the following regularized integral
\begin{equation}
  \label{eq:reg-int}
  \int^{\reg} F(\tau) := \CT_{s = 0}\left[ \lim_{T \to \infty} \int_{\Fc_T} F(\tau) v^{-s} d\mu(\tau)\right].
\end{equation}

Let $k, \ell \in \Nb_0$. 
For a lattice $L$ as in the previous section and
a weakly holomorphic modular form 
$f = \sum_{m, \mu} c(m, \mu) q^m \phi_\mu \in M_{k-\ell+\frac{1}{2},L^-}^!$ 
of weight $k+\frac{1}{2}$ for $\rho_{L^-}$, we define the theta lift
\begin{equation}
  \label{eq:Phikl}
  \Phi^{(\ell, k)}(f, z, h) := \int^\reg \langle f(\tau), \Theta_L^{(\ell, k)}(\tau, z, h)\rangle. 
\end{equation}
When $\ell = 0$, this is the Borcherds-Shimura theta lift
\begin{equation}
  \label{eq:BS}
  \Phi_{\BS}^{}(f,z,h) := \Phi^{(0, k)}(f, z, h)
\end{equation}
studied in \cite{borcherds}. It is a meromorphic modular form of weight $2k$ on $\Gamma_L$ with meromorphic singularities at CM points.
When $f$ has weight $3/2 - k$, there is a related regularized theta lift. This was later studied by Zemel in \cite{Zemel15} and we call it the Borcherds-Zemel lift. 
It is given by
\begin{equation}
  \label{eq:BZtheta}
  \Phi_\BZ(f, z, h) :=
   \Phi^{(2k-1, k)}(f, z, h)
\end{equation}
and defines a meromorphic modular form of weight $2k$
with similar singularities as the Borcherds-Shimura lift.

Let $\varphi_{f, m} \in \Sc_{L^-}$ be the Schwartz function constructed from the $m$-th Fourier coefficients of $f$ in \eqref{eq:varphim}.
There are finitely many $m< 0$ such that $\varphi_{f, m}$ is non-trivial.
Given such an $m < 0$, let $U^- \subset V$ be a negative definite subspace such that $y_U / \sqrt{|m|} \in \Qb^\times$, where $y_U = \Im(z_U)$ with $z_U = z_U^+$ as in \eqref{eq:ZU}.
Then for $z$ is a neighborhood of $z_U$ and $[z_U, t] \in Z(U)$, the modular forms $\Phi_\BS(f, z, t)$ and $\Phi_\BZ(f, z, t)$ have the following Laurent expansions.
\begin{proposition}
  \label{prop:sing}
  In the notations above, for $z$ in a neighborhood of $z_U$, we have
\begin{equation}
  \label{eq:singularity}
  \begin{split}
     (1-w)^{-2k} \Phi^{}_\BS(f, z, t)
&= \frac{2\Gamma(k)}{(4 \sqrt{2|m|} \pi)^k}
\varphi_{f,m}(t^{-1} \lambda_U)
(y_U w)^{-k}
+ O(1),\\
     (1-w)^{-2k} \Phi^{}_\BZ(f, z, t)
&=
\frac{
  \sqrt{|m|}^{k-1} \Gamma(k)}{      \sqrt\pi (\sqrt2)^{k}}
\varphi_{f,m}(t^{-1} \lambda_U)
  (y_U w)^{-k}
+ O(1),
  \end{split}
\end{equation}
where $w = \frac{z - z_U}{z - \overline{z_U}}$ and $\lambda_U = \sqrt{|m|}Z^\perp(z_U) \in V(\Qb)$. 
\end{proposition}
\begin{proof}
  This result follows from specializing Theorem 6.2 in \cite{borcherds}. For completeness, we include a proof here for $\Phi_{\BZ}$.
  It suffices to consider the case $f$ is the Poincar\'e series
  $$
  f_{m, \mu}(\tau) := \sum_{\gamma \in \tilde\Gamma_\infty \backslash \Mp_2(\Zb)} (\ebf(-m \tau) \mid_{3/2-k} \gamma ) \overline{\rho_L}(\gamma)^{-1} \phi_\mu,
  $$
  where $\tilde\Gamma_\infty \subset \Mp_2(\Zb)$ is the preimage of $\Gamma_\infty = \{\smat{1}{b}{}{1}: b \in \Zb\} \subset \SL_2(\Zb)$ under the natural surjection. 
  Unfolding the integral against it, we obtain
  \begin{align*}
    \Phi_\BZ&(f, z, t)
    = \int^{\mathrm{reg}}_{\Gamma_\infty \backslash \Hb}
      \ebf(-m \tau) \langle \phi_\mu + (-1)^{k+1}\phi_{-\mu}, \Theta_L^{(2k-1, k)}(\tau, z ,t)\rangle d\mu(\tau)\\
    &=
\frac{      2^{1-2k}}{      (-\sqrt2 y^2)^{k}}
      \lim_{T \to \infty}
       \int_{0}^T
    v^{3/2}
    \sum_{\substack{\lambda \in t(L\pm\mu)\\ Q(\lambda) = m}}
(\pm 1)^{k+1} \He_{2k-1} \lp \sqrt{\pi v}
    \lp \lambda, Z^\perp(z)\rp\rp
  \left(\lambda,
    Z(\bar z)
    \right)^{k}
e^{4\pi vQ(\lambda_z)}
    \frac{dv}{v^2}\\
    &=
      \frac{(\phi_\mu(t^{-1} \lambda_U) +   (-1)^{k+1}\phi_{-\mu}(t^{-1} \lambda_U) + \phi_\mu(-t^{-1} \lambda_U) +   (-1)^{k+1}\phi_{-\mu}(-t^{-1} \lambda_U) )}{      (-\sqrt2 y^2)^{k}}\\
&\quad \times 2^{1-2k}  \left(\lambda_U,
    Z(\bar z)
    \right)^{k}     \int_{0}^\infty
 \He_{2k-1} \lp \sqrt{\pi v}
    \lp \lambda_U, Z^\perp(z)\rp\rp
e^{-\frac{ \pi v}{y^2}|(\lambda_U, Z(z))|^2}
\sqrt{v}      \frac{dv}{v}
       + O(1)
  \end{align*}
  for $z \neq z_U$ in a neighborhood of $z_U$ with $[z_U, t] \in Z(U)$.
  Now for any $1 \le j \le k$, we have
  \begin{align*}
y^{-2k}  \left(\lambda_U,
    Z(\bar z)
    \right)^{k}
&    \lp \lambda_U, Z^\perp(z)\rp^{2j-1}
    \int_{0}^\infty
    v^{j}
e^{-\frac{\pi v}{y^2}|(\lambda_U, Z(z))|^2}
  \frac{dv}{v}
  = (\bar w)^{k-j} \cdot P_j(w, \bar w),
  \end{align*}
  with $P_j \in \Cb(\!(w)\!) \llbracket \bar w \rrbracket$ a power series with coefficients in Laurent series in $w$.
  Furthermore applying \eqref{eq:change-var} and $\lambda_U = \sqrt{|m|} Z^\perp(z_U)$,   it is straightforward to evaluate
  \begin{align*}
P_k(w, \bar w)
&= 
  \frac{\sqrt{|m|}^{k-1}\Gamma(k)}{\pi^k}
  \frac{(Z^\perp(z_U), Z^\perp(z))^{2k-1}}{  (Z^\perp(z_U), Z(z))^{k}}
  =  \frac{\sqrt{|m|}^{k-1}\Gamma(k)}{2(-\pi)^k}
\frac{(1-w)^{2k}}{  (y_U w)^{k} }
+ \bar w R(w, \bar w)
  \end{align*}
for some $R(w, \bar w) \in   \Cb(\!(w)\!) \llbracket \bar w \rrbracket$.
Therefore, we can apply \eqref{eq:parity} to obtain
$$
\Phi_\BZ(f, z, t)
-
\frac{(\phi_\mu(t^{-1} \lambda_U) + (-1)^{k-1}\phi_\mu(-t^{-1} \lambda_U))
  \sqrt{|m|}^{k-1} \Gamma(k)}{      \sqrt\pi (\sqrt2)^{k}}
\frac{(1-w)^{2k}}{ (y_U w)^{k}} = \bar w P(w, \bar w) + O(1)
$$
for some $P(w, \bar w) \in   \Cb(\!(w)\!) \llbracket \bar w \rrbracket$.
  Since $\Phi_\BZ(f, z, t)$ is meromorphic in $z$ for all $t$, we can set $\bar w = 0$ in the equation above gives us \eqref{eq:singularity}. 
\end{proof}

To describe the Laurent expansion of $\Phi_\BS$, we need the following result that relates the constant term of the Laurent expansion of $\Phi^{(\ell, k)}$ and its ``value''. This is an example of the over-regularization phenomenon for regularized theta lifts (see \cite[section 1]{Schofer}).

\begin{proposition}
  \label{prop:over-reg}
  Let $\ell \in \Nb_0$ and $k \in \Nb$.
  For any $(z_0, h_0) \in X_K$, the regularized integral defining $\Phi^{(\ell,k)}$ exists at $(z_0, h_0)$ and equals to $\CT^{(2k)}_{z_0}(\Phi^{(\ell , k)})$. 
\end{proposition}

\begin{proof}
  We only need to verify the claim when $[z_0, h_0] = [z_U, t] \in Z(U)$ is on the singularity of $\Phi^{(\ell, k)}$.
  Without loss of generality we take $t = 1$ and omit it from the notation.
  Let $w =\frac{z-z_U}{z-\overline{z_U}}$ be the local coordinate as in \eqref{Awz}.
    When $|w|$ is small but nonzero, we have
    \begin{align*}
  \Phi^{(\ell, k)}(f, z)
&=
\CT_{s = 0}
\lp
\lim_{T \to \infty} \int_{1}^T
\sum_{\mu \in L'/L} \CT(f_\mu(\tau) \theta_V^{(\ell, k)}(\tau, z, \phi_\mu)) v^{-s}\frac{dv}{v^2}\rp\\
&+ \int_{\Fc_1} \langle f(\tau), \Theta_L^{(\ell ,k)}(\tau, z)\rangle d \mu(\tau).
    \end{align*}
    The second term exists for any $z$. 
Up to the constant $2^\ell(-\sqrt{2})^k$, the limit in the first term above becomes
  $$
 y^{-2k}
  \lim_{T \to \infty} \int_{1}^T
\sum_{\substack{m < 0\\\lambda \in V,~ Q(\lambda) = -m}} \varphi_m(\lambda)
 \He_\ell(\sqrt{\pi v}(\lambda, Z^\perp(z))
(\lambda, Z(\bar z))^k e^{-4\pi v |Q(\lambda_z)|}
v^{k - \ell/2-s}\frac{dv}{v}.
$$
The contribution to limit above from  $\lambda \neq \pm\sqrt{m}Z^\perp(z_U)$ exists for any $z$ in a small neighborhood of $z_U$, including $z = z_U$.
Using \eqref{eq:change-var}, the contribution
of $\lambda = \pm\sqrt{m}Z^\perp(z_U)$ to the limit is
$$
(\pm 1)^{k + \ell}\sqrt{m}^k \varphi_m(\lambda)
 y^{-2k}
  \lim_{T \to \infty} \int_{1}^T
 \He_\ell\lp 2 \sqrt{m\pi v} \frac{1 + |w|^2}{1-|w|^2}\rp
\lp \frac{-4y_U \bar w}{ ( 1 - \bar w)^2} \rp^k e^{-\epsilon v}
v^{k - \ell/2-s}\frac{dv}{v}
$$
with $\epsilon = 4\pi|Q(\lambda_z)| = 16m \pi |w|^2/ (1 - |w|^2)^2$. 
When $w =0$, the integrand is identically zero since $k \ge 1$.
This shows that the regularized integral defining $\Phi^{(\ell, k)}$ exists at $z = z_U$.
When $w \neq 0$, the integral is well-defined for $s = 0$ and we have
$$
\CT\lp
(1-w)^{-2k}
 y^{-2k}
  \lim_{T \to \infty} \int_{1}^T
 \He_\ell\lp 2 \sqrt{m\pi v} \frac{1 + |w|^2}{1-|w|^2}\rp
\lp \frac{-4y_U \bar w}{ ( 1 - \bar w)^2} \rp^k e^{-\epsilon v}
v^{k - \ell/2}\frac{dv}{v}\rp = 0.
$$
This proves that $\Phi^{(\ell, k)}(f, z_U) = \CT^{(2k)}_{z_U}(\Phi^{(\ell, k)}(f, z))$. 
\end{proof}

\subsection{Differential operators on theta functions}
In this section we study the action of the raising operator on the theta function $\theta_V^{(\ell,k)}$. To this end, we recall some differential operator calculations in the Fock model of the Weil representation (see e.g.\ \cite[Appendix]{Li22-average}).
Let  $\varphi^\circ \in \Sc(\Rb^{1, 2})$ be the Gaussian
$$
\varphi^\circ(x_1, x_2, x_3) := e^{-\pi(x_1^2 + x_2^2 + x_3^2)}.
$$
and $D_r$ are operators on $\Sc(\Rb^{1, 2})$ defined by
\begin{equation}
  \label{eq:Dj}
D_r := \partial_{x_r} - 2\pi x_r,~ 1 \le r \le 3.
\end{equation}
Let $\mathbb{S}(\Rb^{1, 2}) \subset \Sc(\Rb^{1, 2})$ be the subspace spanned by  functions of the form $\prod_{1 \le j \le 3} D_j^{r_j} \varphi^\circ$ for $r_j \in \Nb_0$. 
It is called the polynomial Fock space and isomorphic to $\Cb[\zf_1, \zf_2, \zf_3]$ via
\begin{equation}
  \label{eq:iota}
  \iota: \mathbb{S}(\Rb^{1, 2}) \to \Cb[\zf_1, \zf_2, \zf_3],~
  D_1^{r_1}   D_2^{r_2}   D_3^{r_3} \varphi^\circ \mapsto i^{r_1 - r_2 - r_3} \zf_1^{r_1}
  \zf_2^{r_2}\zf_3^{r_3}.
\end{equation}
 We now set
 \begin{equation}
   \label{eq:vwf}
   \wf := \zf_2 - i\zf_3 .
 \end{equation}
 It is easy to check that the Schwartz function $\varphi^{(\ell, k)}_\infty$ defined in \eqref{eq:varphik} is
 \begin{equation}
   \label{eq:varphilkfock}
 \varphi^{(\ell, k)}(x_1, x_2, x_3) =
(-2\sqrt{2\pi})^{-\ell}
 (-4\pi)^{-k} 
 D_1^\ell  (D_2 - i D_3)^k \varphi^\circ(x_1, x_2, x_3),
 \end{equation}
 and is in $\mathbb{S}(\Rb^{1, 2})$, and its image under $\iota$ is given by
 \begin{equation}
   \label{eq:varphilk}
   \iota(\varphi^{(\ell, k)}_\infty) =
(2\sqrt{2\pi} i)^{-\ell}
 (-4\pi i )^{-k} 
  \zf_1^\ell \wf^k. 
 \end{equation}
 The elements $L, R \in \slf_2(\Cb) \cong \spf(W \otimes \Cb)$ corresponding to the lowering and raising operator in the symplectic variable
 act on $\Cb[\zf_1, \zf_2, \zf_3]$ as (see \cite[Lemma A.2]{FM06})
\begin{equation}
  \label{eq:fock1}
      \omega(L) = - 2\pi  \partial_{\zf_1}^2 +    \frac{1}{8\pi} \wf \overline{\wf},~
  \omega(R) = - 8\pi \partial_{\wf} \partial_{\overline{\wf}} + \frac{1}{8\pi} \zf^2_1.
\end{equation}
The images of the elements $L, R \in   \slf_2(\Cb) \cong \mathfrak{so}_3(\Cb)$, which we denote by $\Lz, \Rz$, act on $\Cb[\zf_1, \zf_2, \zf_3]$ through $\iota$ as (see \cite[Lemma A.1]{FM06})
\begin{equation}
  \label{eq:fock}
  \begin{split}
    \omega(\Lz) &= 8\pi  \partial_{\zf_1} \partial_{{\wf}} -    \frac{1}{4\pi} {\zf_1} \overline{\wf},~
  \omega(\Rz) =  8\pi \partial_{{\zf_1}} \partial_{\overline{\wf}} - \frac{1}{4\pi} \zf_1 \wf,
  \end{split}
\end{equation}
For convenience, we slightly abuse notation and write $L, R, \Lz, \Rz$ for their corresponding actions on $\Cb[\zf_1, \zf_2, \zf_3]$.
Then for $k, \ell\ge 0$, we have
\begin{equation}
  \label{eq:diff-op-varphi}
  R\varphi^{(\ell, k)} = -\varphi^{(\ell+2, k)},~
  R_z\varphi^{(\ell, k)} = 2\sqrt{2\pi }  \varphi^{(\ell+1, k+1)}.
\end{equation}
In terms of theta functions, we have
\begin{equation}
  \label{eq:diff-op-theta}
  \begin{split}
      R_\tau \theta_V^{(\ell, k)}
 & =
   -      \theta_V^{(\ell+2, k)},~
   R_z \theta_V^{(\ell, k)}
 =
2\sqrt{2\pi }        \theta_V^{(\ell+1, k+1)}.
  \end{split}
\end{equation}
    From this, we have the following result.

\begin{proposition}
  \label{prop:diff-id}
  For any $m, k \in \Nb$, we have
  \begin{equation}
    \label{eq:raise-theta}
    \begin{split}
          R_{z}^m \theta_V^{(0, k)}
&    =
    (    2\sqrt{2\pi } )^{m} \theta_V^{(m, k+m)}
    =
(    2\sqrt{2\pi } )^{m} (-1)^r
R_{\tau}^r \theta_V^{(\ell, k+m)},\\
          R_{z}^m \theta_V^{(2k-1, k)}
&    =
    (    2\sqrt{2\pi } )^{m} \theta_V^{(2k-1+m, k+m)}
    =
(    2\sqrt{2\pi } )^{m} (-1)^{k+r+\ell -1}
R_{\tau}^{k+r+\ell-1} \theta_V^{(1-\ell, k+m)}
    \end{split}
  \end{equation}
  where $m = 2r + \ell$ with $r \in \Nb$ and $\ell \in \{0, 1\}$. 
\end{proposition}

\section{Recurrences for Laurent coefficients of modular forms}
\label{sec:rec}
When~$f$ is a holomorphic modular form, formula~\eqref{eq:LaurentC} says that~$a_m=\frac{y_0^m\bigl(R_k^mf\bigr)(z_0)}{m!}$. 
If~$z_0$ is a CM point of discriminant $-D< 0$ and $f$ has algebraic Fourier coefficients, it is known that~$\bigl(R_k^mf\bigr)(z_0)$ is of the form~$\Omega_{-D}^{k+2m}\cdot\alpha_m(f)$, where~$\Omega_{-D}$ is the Chowla-Selberg period defined in \eqref{eq:Chowla-Selberg}, and~$\alpha_m(f)$ is an algebraic number that depends on~$f$. 
In~\cite{villegaszagier} 
Rodriguez--Villegas and Zagier showed that the numbers~$\{\alpha_m(f)\}_m$ are special values at an algebraic point of a family of polynomials~$\{p_m(t)\}_m$ that satisfy a\;\lq\lq{quasi-recursion}\rq\rq. In general, the quasi-recursion looks like
\begin{equation}
\label{eq:grec}
(m+1)p_{m+1}(t)=a_1(m,k,t)\frac{dp_m(t)}{dt}+a_2(m,k,t)p_m(t)+a_3(m-1,k,t)p_{m-1}(t)\,,
\end{equation}
for some polynomials~$a_i(x,y,t)\in\Qb[x,y,t]$ that we specify below (this is called \emph{quasi-recursion} because of the presence of the derivative of~$p_m(t)$).

The main idea behind the quasi-recursive method is the following. Fix a CM point~$z_0$ and let~$f\in M_k(\SL_2(\Zb))$. Let~$\varphi$ be a (possibly meromorphic) quasimodular form of weight two such that~$\varphi^*(z):=\varphi(z)-\frac{1}{4\pi y}$, where~$y=\mathrm{Im}(z)$, has a zero in~$z_0$.
Consider the~\emph{Serre derivative}, defined recursively by
\begin{equation}
\label{eq:Serre}
\vartheta_\varphi^{[0]}{f}=f\,,\quad\vartheta_\varphi^{[1]}=\vartheta_\varphi{f}:=\frac{1}{2\pi i}\partial_z{f}-k\varphi{f}\,,\quad
\vartheta_\varphi^{[m+1]}=\vartheta_\varphi\bigl(\vartheta_\varphi^{[m]}{f}\bigr)-n(n+k-1)\Phi\vartheta_\varphi^{[m-1]}{f}\,,
\end{equation}
where~$\Phi:=\varphi'-\varphi^2$ is a modular form of weight~$4$. The function~$\vartheta_\varphi^{[m]}f$ is a holomorphic modular form of weight~$k+2m$. 
The relation between the Serre derivative~$\vartheta_\varphi$ and the raising operator is
\begin{equation}
\label{eq:RaisingSerre}
R_k^m{f}(z)=(-4\pi)^m\sum_{r=0}^m{\binom{m}{r}(k+r)_{m-r}\bigl(\varphi^*(z)\bigr)^{m-r}\vartheta_\varphi^{[r]}f(z)}\,,
\end{equation}
Since~$\varphi^*(z_0)=0$, the identity~\eqref{eq:RaisingSerre} implies that~$R_k^mf(z_0)=(-4\pi)^m\vartheta_\varphi^{[m]}f(z_0)$, so the computation of the coefficients~$\alpha_m(f)$ is reduced to the computation of~$\vartheta_\varphi^{[m]}f(z_0)$.

Suppose that~$z_0\neq e^{2\pi i/3}$ and set~$t:=E_6E_4^{-3/2}$ (for ~$z_0=e^{2\pi i/3}$ we can take $t = E_4/E_6^{2/3}$).
Let~$F_m$ be the polynomial such that~$\vartheta_{\varphi}^{[m]}f=F_m(E_4,E_6)$. Then~$F_m(E_4,E_6)E_4^{-(k+2m)/4}$ has weight zero and is by construction a polynomial in~$t$ that we denote~$p_m(t)$. Note that, if~$t_0=t(z_0)$, 
\begin{equation}
\label{eq:normO}
\vartheta_\varphi^{[m]}f(z_0)=p_m(t_0)\cdot E_4(z_0)^{\frac{k+2m}{4}}=p_m(t_0)\cdot(c_0\Omega_{-D})^{k+2m}
\end{equation}
for some computable algebraic number~$c_0$
The quasi-recursion~\eqref{eq:qrec} for the polynomials~$p_m(t)$ eventually follows from the recursive definition of the Serre derivative~\eqref{eq:Serre} as follows and the definition of the~$p_m(t)$.
For details see Sections~6 and 7 of~\cite{villegaszagier} or Section 6.3 of the first part of~\cite{brdgza08}.  
One deduces from the above construction that the polynomial coefficients~$a_i(x,y,t)$ of the recursion~\eqref{eq:qrec} are of the form
\begin{equation}
\label{eq:ai}
a_1(x,y,t)=a_1(t)\,,\qquad a_i(x,y,t)\;=\;x\cdot a_{i,1}(t)\;+\; y\cdot a_{i,2}(t)\,,\quad i=2,3\,.
\end{equation}
As an example, when~$z_0=i$ one can choose~$\varphi=E_2$. It follows that~$\Phi=-E_4/144$ and that the polynomials~$\{p_m(t)\}_m$ satisfy the quasi-recursion
\begin{equation}
\label{eq:qrec}
(m+1)p_{m+1}(t)=\frac{(t^2-1)}{2}\frac{dp_m(t)}{dt}-\frac{(k+2m)t}{12}p_m(t)-\frac{(m-1+k)}{144}p_{m-1}(t)\,.
\end{equation}

Now ~$g$ be a meromorphic modular form with algebraic Fourier coefficients and  a pole of order~$N$ at~$z_0$. The Laurent coefficients of~$g$ in~$z_0$ are still expected to be products of a period and an algebraic number, but the argument above does not work to give a quasi-recursive family of polynomials that specialize to the algebraic part of the coefficients. This is because the constant term of~$R_k^m{g}$ is not related to the value or to the constant term of the Serre derivative~$\vartheta^{[m]}g$. This can be seen from equation~\eqref{eq:RaisingSerre} (with~$f$ replaced by~$g$): the order of the pole of~$\vartheta^{[r]}{g}$ in~$z_0$ is~$N+r$, while~$(\varphi^*)^{m-r}$  has a zero of order~$m-r$ there; it follows that the constant term on the right hand-side of~\eqref{eq:RaisingSerre} will not in general be the constant term of the expansion of~$\vartheta_\phi^{[m]}{g}$, but some linear combination involving also Serre derivatives of lower degree and powers of~$\varphi^*$. For such a combination there is no modular interpretation that leads to a simple recursive formula. 

In the following we construct two different families of polynomials (or rational functions) whose evaluation at an algebraic point give the algebraic part of Laurent coefficients of~$g$. We present them in the case of~$\SL_2(\Zb)$ but, as for the Rodriguez--Villegas--Zagier method, the construction extends to other arithmetic groups.
The main tool is the interpretation of the recursion~\eqref{eq:grec} in terms of the linear differential operator
\begin{equation}
\label{eq:L}
\L_k\;=\;a_1(t)\frac{d}{dt}\;+\;\bigl(a_{3,1}(t)X^2+a_{2,1}(t)X-1\bigr)\frac{d}{dX}\;+\;\bigl(a_{3,2}(t)X+a_{2,2}(t)\bigr)k\,,
\end{equation}
where~$a_{i,j}(t)$, for~$i=2,3$, are defined in~\eqref{eq:ai}. 
A straightforward computation gives us the following lemma.
\begin{lemma}
\label{lem:L}
It holds~$\L_k\bigl(\sum_{m\ge0}{p_m(t)X^m}\bigr)=0$ if and only if~$\{p_m(t)\}_{m\ge0}$ satisfy the recursion~\eqref{eq:grec}. 
\end{lemma}


The first method is given in the next proposition.
\begin{proposition}
\label{prop:rec1}
Let~$z_0=x_0+iy_0$ be a CM point with associated Chowla-Selberg period~$\Omega_{-D}$. 
Let~$g$ be a meromorphic modular form of weight~$k \in 2\Zb$ on~$\S\L_2(\Zb)$ with a pole of order~$N\ge1$ at~$z_0$. Let~$a_i(x,y,t)\,,i=1,2,3$ be as in~\eqref{eq:grec}.
Then the Laurent expansion of~$g$ at~$z_0$ is given by
\[
(1-w)^{-k}g\biggl(\frac{w\bar{z}_0-z_0}{w-1}\biggr)=\sum_{m=-N}^\infty{q_m(t_0)\Omega^{k+2m}(-4\pi y_0 w)^m}
\]
where~$\Omega$ is an known algebraic multiple of~$\Omega_D$, $t_0\in\bar{\Qb}$, and the rational functions~$q_m(t)$ satisfy the following non-linear recursion
\begin{align*}
(m+1+2N)q_{m+1}(t)&=a_1(t)q_n'(t)+a_2(m,k,t)q_m(t)+a_3(m-1,k,t)q_{m-1}(t)\\
&-A(t)\sum_{j=-N}^{m}{q_j(t)q_{m-N-j}(t)}-B(t)\sum_{j=1-N}^m{q_j(t)q_{m+1-N-j}(t)}\,.
\end{align*}
with
\[
A(t)=\frac{(N+1)}{q_{-N}(t)^2}\bigl(q_{1-N}(t)+a_1(t)q_{-n}'(t)-a_2(N,-k,t)q_{-N}(t)\bigr)\,,\quad B(t)=\frac{N}{q_{-N}(t)}\,,
\]
and initial data dependent on~$g$.
\end{proposition}
\begin{remark}
The initial data~$q_{-N}$ and~$q_{1-N}$ are explicitly given in the proof in \eqref{eq:AB}.
\end{remark}
\begin{proof}
The function~$g^{-1}$ is a modular form of weight~$-k$ holomorphic in~$z_0$ and possibly meromorphic at other points of~$\Hb\cup\infty$. Let~$\{p_m(t)\}$ be the sequence of polynomials constructed from the recursion~\eqref{eq:grec} with any initial data such that~$p_N(t)\neq0$. By Lemma~\ref{lem:L}, the generating series~$P(t,X):=\sum_{m\ge0}{p_m(t)}X^m$ satisfies~$\L_{-k}P(t,X)=0$.
Now consider the series~$Q(t,X)=\sum_{m=-N}^\infty{q_m(t)X^m}$ defined by
\[
\frac{1}{Q(t,X)}:=P(t,X)-\sum_{m=0}^{N-1}{p_m(t)X^m}\,.
\]
For every~$m\ge -N$, the coefficients of~$q_m(t)$ of~$Q(t,X)$ are rational functions in~$t$ holomorphic in~$t=0$. This is because
\begin{equation}
\label{eq:Q}
Q(t,x)=\frac{1}{p_N(t)X^N\bigl(1+\sum_{m=1}^\infty{\frac{p_{N+m}(t)}{p_N(t)}X^m}\bigl)}
\end{equation}
and~$p_N(0)\neq0$ by hypothesis. 
The computation
\begin{align*}
2\bigl(a_{3,2}(t)X+a_{2,2}(t)\bigr)k&=\bigl(a_{3,2}(t)X+a_{2,2}(t)\bigr)k-\L_{-k}(1)\\
&=\bigl(a_{3,2}(t)X+a_{2,2}(t)\bigr)k-\L_{-k}\Bigl(Q(t,X)\Bigl(P(t,X)-\sum_{m=0}^{N-1}p_m(t)X^m\Bigr)\Bigr)\\
&=\L_{-k}\Bigl(\sum_{m=0}^{N-1}p_m(t)X^m\Bigr)Q(t,X)-\L_{-k}Q(t,X)\Bigl(P(t,X)-\sum_{m=0}^{N-1}p_m(t)X^m\Bigr)
\end{align*}
together with the identity
\[
\L_{-k}\Biggl(\sum_{m=0}^{N-1}p_m(t)X^m\Biggr)=Np_N(t)X^{N-1}+a_3(N-1,-k,t)p_{N-1}(t)X^N\,,
\]
which follows from~\eqref{eq:grec} and the definition of~$\L_k$, 
shows that the series~$Q(t,X)$ satisfies the non-linear differential equation
\begin{equation}
\label{eq:diff1}
\L_kQ(t,X)\;-\;\Bigl(Np_N(t)X^{N-1}+a_3(N-1,-k,t)p_{N-1}(t)X^N\Bigl)Q(t,X)^2=0\,.
\end{equation}
Since, by definition of~$Q(t,X)$, 
\[
p_N(t)=q_{-N}(t)^{-1}\,,\qquad p_{N-1}(t)=\frac{(N+1)}{q_{-N}(t)^2}\bigl(q_{1-N}(t)+a_1(t)q_{-N}'(t)-a_2(N,-k,t)q_{-N}(t)\bigr)\,,
\]
the above differential equation translates in the non-linear recursion for the coefficients~$q_n(t)$ of~$Q(t,X)$ given in the statement.

Finally, as initial values for the recurrence~\eqref{eq:grec} choose~$p_{-1}(t)=0$ and~$p_0(t)$ the polynomial associated to the holomorphic modular form~$g^{-1}$ by the Rodriguez--Villegas--Zagier method. There exists a~$t_0\in\bar{\Qb}$ (more precisely, for a fixed branch one has $t_0=t(z_0)=E_6(z_0)E_4(z_0)^{3/2}$ if~$z_0\neq e^{2\pi i/3}$ or~$t_0=E_4(z_0)E_6(z_0)^{2/3}$ otherwise) such that~$P(t_0,-4\pi y_0\Omega^2w)\Omega^{-k}$ is the Laurent series of~$g^{-1}$ in~$z_0$, where~$\Omega=c_0\cdot\Omega_{-D}$ for an algebraic number~$c_0$ (see~\eqref{eq:normO}).
If 
\begin{equation}
\label{eq:AB}
q_{-N}(t)=\frac{1}{p_N(t)}\,,\quad q_{1-N}(t)=\frac{p_{N+1}(t)}{p_N(t)^2}\,,
\end{equation}
then~$Q(t_0,-4\pi y_0\Omega^2w)\Omega_{-D}^k$ is the Laurent expansion of~$g$ because
\[
Q\bigl(0,-4\pi y_0\Omega_{-D}^2w\bigr)=P\bigl(0,-4\pi y_0\Omega_{-D}^2w\bigr)^{-1}
\] 
and~$p_n(t_0)=0$ for~$n=0,\dots,N-1$ since~$g^{-1}$ has a zero of order~$N$ in~$z_0$. 
 
\end{proof}

\begin{remark}
\label{rmk:rec1}
\begin{enumerate}[wide=0pt]
\item If one defines~$Q(t,X)^{-1}=P(t,X)-\sum_{m=0}^N{p}_m(t)X^m+{p}_N(0)$ in the proof of Proposition~\ref{prop:rec1}, then~$Q(t,X)$ has polynomial coefficients in~$t$ (if~$p_n(t)$ for~$n\ge N$ are polynomials), and~$Q(t,X)$ satisfies a slightly more complicated recursion formula.
\item In the case~$N=1$, which turns out to be of relevance in view of Theorem~\ref{thm:rec2}, the polynomials~$q_{-1}(t)$ and~$q_0(t)$ are easily determined from~$g$, similarly to the holomorphic case. Thanks to~\eqref{eq:diff1}, the rational functions~$A(t)$ and~$B(t)$ in the statement of~Proposition~\ref{prop:rec1} can be written as~$A(t)=a_3(0,-k,t)p_0(t)$ and~$B(t)=p_1(t)$, where~$p_0(t)$ and~$p_1(t)$ are the polynomials obtained as in the Rodriguez--Villegas--Zagier method from the holomorphic modular forms~$g^{-1}$ and~$\vartheta_\varphi(g^{-1})$ respectively. More generally, in the computations it is more efficient to use the description~$A(t)=a_3(N-1,-k,t)p_{N-1}(t)$ and~$B(t)=Np_N(t)$ deduced from~\eqref{eq:diff1} than the one in the statement of~Proposition~\ref{prop:rec1}, since one has to compute~$p_n(t)$ up to~$n=N$ instead of~$n=N+1$.
\end{enumerate}
\end{remark}

\begin{example}
\label{ex:rec1}
We compute the Laurent expansion of~$-2^8\frac{\Delta}{E_4^2}$ in~$\zeta = \frac{-1+ \sqrt3 i}{2}$, where we have~$k=4$ and~$N=2$ (a double pole in~$\zeta$). 
For the CM point~$\zeta$ one can still choose~$\varphi=E_2$ and the polynomials~$a_i(x,y,t)$ in~\eqref{eq:grec}, if~$t=E_4E_6^{-2/3}$ are
\[
a_1(t)=\frac{t^3-1}{2}\,,\quad a_2(x,y,t)=-\frac{(2x+y)t^2}{12}\,,\quad a_3(x,y,t)=-\frac{(x+y)t}{144}\,.
\]

The recurrence in Proposition~\ref{prop:rec1} is determined by the rational functions $A(t)=-\frac{t p_1(t)}{48}$ and~$B(t)=2p_2(t)$ where~$p_1(t)$ and~$p_2(t)$ are obtained from the recursion~\eqref{eq:qrec} with initial data~$p_{-1}(t)=0$ and~$-2^{-8}\frac{E_4^2}{\Delta}=\frac{27t^2}{4(t^3-1)}=p_0(t).$ It follows that~$A(t)=\frac{3t^2}{32(t^3-1)}$ and~$B(t)=\frac{39t^2+24}{32(t^3-1)}$. Normalizing as in Remark~\ref{rmk:rec1}, the first terms we get from the non-linear recursion are~$q_{-2}(t)=-4/3\,,q_{-1}(t)=\frac{11t^2}{9(t^3 -1)}\,,q_0(t)=\frac{-277t^7 - 2019t^4 + 360t}{1728(t^3-1)^2}\,,q_1(t)=\frac{19699t^9 + 49597t^6 - 16776t^3 + 720}{51840(t^3-1)^3}$.
Since~$t(\zeta)=t_0=0$, and taking into account the identity ~$E_6(\zeta)=24\sqrt{3}\Omega_{-3}^6$, Proposition~\eqref{prop:rec1} gives the same Laurent expansion as in~\eqref{eq:ex1}. 
\end{example}

In contrast with Proposition~\ref{prop:rec1}, Theorem~\ref{thm:rec2} exhibits a family of polynomials~$\{q_m(t)\}_m$ that satisfy a \emph{linear} recursion and whose evaluation at an algebraic point gives the Laurent coefficients of~$g$. In order to compute it, one needs the Laurent expansion of the logarithmic derivative of a suitable modular function~$h$; since by construction it has only a simple pole, such expansion can be easily computed with Proposition~\ref{prop:rec1}. The choice of~$h$ depends only on the point~$z_0\in\Hb$ and is universal: the Laurent expansion of~$\tfrac{d}{d\tau}\log(h)$ has to be computed once for all, and can be used to find the Laurent expansion of every modular form with a pole in~$z_0$ (of any order).  Again, we give the complete result in the case of~$\S\L_2(\Zb)$, but the same idea applies to other arithmetic groups.

\begin{proof}[Proof of Theorem~\ref{thm:rec2}]
If~$N=0$ then~$g$ is holomorphic in~$z_0$ and we are done. Assume than that~$N\ge1$ and let~$h$ be a modular function (locally around~$z_0$) with a simple zero in~$z_0$.
The product~$f=gh^N$ is a modular form of weight~$k$ holomorphic in~$z_0$. Let~$\{p_m(t)\}$ be a sequence of polynomials constructed from the recursion~\eqref{eq:qrec} with any initial data, and consider the generating series~$P(t,X):=\sum_{m\ge0}{p_m(t)}X^m$. 
Let~$\{\tilde{h}_m(t)\}_m$ be the sequence of polynomials obtained from the recursion~\eqref{eq:qrec} with initial data depending on~$h$ (as in the Rodriguez--Villegas--Zagier method) and consider the generating function~$H(X)=\sum_{m\ge1}{\tilde{h}_m(0)X^m}$. 

Define~$Q(t,X):=P(t,X)\cdot H(X)^{-N}$.
The coefficients~$q_m(t)$ of~$Q(t,X)=\sum_{m\ge -N}{q_m(t)X^m}$ are polynomials if the~$\{p_m(t)\}_m$ are.
From Lemma~\ref{lem:L} one has
\begin{equation}
\label{eq:rec2}
0=\L_k P(t,x)=\L_k\bigl(Q(t,X)H(X)^N\bigr)=\L_k\bigl(Q(t,X)\bigr)H(X)^N+\L_0\bigl(H(X)^N\bigr)Q(t,X)\,.
\end{equation}
Since~$H(X)$ does not depend on~$t$, the expression~$\L_0(H(X)^N)$ reduces to~$(a_{3,1}(t)X^2+a_{2,1}(t)X-1)\frac{d}{dX}H(X)^N$. Together with~\eqref{eq:rec2}, this implies that
\begin{align*}
0&=\L_k\bigl(Q(t,x)\bigr)+Q(t,x)\bigl(a_{3,1}(t)X^2+a_{2,1}(t)X-1\bigr)\frac{1}{H(x)^N}\frac{dH(x)^N}{dx}\\
&=\L_k\bigl(Q(t,x)\bigr)+N\cdot Q(t,x)\bigl(a_{3,1}(t)X^2+a_{2,1}(t)X-1\bigr)\frac{d\log(H(X))}{dX}\\
&=\L_k\bigl(Q(t,x)\bigr)+N\cdot Q(t,x)\bigl(a_{3,1}(t)X^2+a_{2,1}(t)X-1\bigr)\sum_{m\ge-1}{(-1)^mh_mX^m}\,,
\end{align*}
from which the recursion in the statement follows (note that~$h_{-1}=1$). 

The solution of the recursion is completely determined if we know~$q_{-N}(t)$ and assume~$q_{-1-N}(t)=0$; from the definition of~$Q(t,X)$ we see that~$q_{-N}=p_0(t)\tilde{h}_1^{-N}$. Moreover, similarly to the proof of~Proposition~\ref{prop:rec1}, the Laurent expansion of~$g$ in~$z_0$ is given by~$Q(t_0,-4\pi y_0 \Omega^2w)\Omega^k$ if~$p_0(t)$ is chosen depending on~$f$ and~$t_0$ and~$\Omega$ are as in the paragraph above~\eqref{eq:AB}. In particular, $p_0(t)$ is the expression in~$t=E_6E_4^{-3/2}$ (if~$z_0\neq e^{2\pi i/3}$, otherwise~$t=E_4E_6^{-2/3}$) of the polynomial in~$E_4,E_6$ associated to~$f=gh^N$ multplied by~$E_4^{-k/4}$ (or $E_6^{k/6}$ respectively).
\end{proof}

\begin{example}
Consider again the modular form~$g=-2^8\frac{\Delta}{E_4^2}$ in~$\zeta = \frac{-1+ \sqrt3 i}{2}$, where we have~$k=4$ and~$N=2$ (a double pole in~$\zeta$). 
To compute with the recurrence in Theorem~\ref{thm:rec2} one needs the Laurent expansion of~$d\log(h)$ and one can choose~$h=t=E_4E_6^{-3/2}$. Then~$d\log(h)=\frac{E_4^2}{3E_6}-\frac{E_6}{3E_4}$. The Laurent expansion of~$d\log(h)$ can be easily computed with Proposition~\ref{prop:rec1} since it has a simple pole: the polynomials defining the non-linear recursion are easily seen to be~$A(t)=\frac{-6t^2}{(t^3-1)}$ and~$B(t)=\frac{-(3t^3+2)}{2(t^3-1)}$, and the expansion of~$d\log(h)$ is as in the statement of Theorem~\ref{thm:rec2} with~$h_{-1}=1,h_0=0,h_1=0,h_2=41/1728,\dots$. One has also compute the coefficient~$\tilde{h}_1=-1/3$ of the Laurent expansion of~$h$ as explained in the last part of the proof of~Theorem~\ref{thm:rec2}.
We remark that this step has to be performed only once, as the expansion of~$d\log(h)$ can (and should) be used to compute the Laurent expansion of all modular forms in~$\zeta$. 
The first term~$q_{-2}(t)$ of the expansion of~$g=-2^8\frac{\Delta}{E_4^2}$ is given by~$q_{-2}(t)=(-\frac{1}{3})^{-2}p_0(t)$, where~$p_0(t)$ is the expression in~$t=E_4E_6^{-2/3}$ of~$gh^2=-2^8\Delta E_6^{-3}$. It follows that~$q_{-2}(t)=9\cdot\frac{4}{27}(t^3-1)$. Finally, the quasi-recurrence in Theorem~\ref{thm:rec2} gives~$q_{-1}(t)=-\frac{8}{9}t^2(t^3-1)\,,q_0(t)=-\frac{14}{27}t^7 + \frac{5}{6}t^4 - \frac{17}{54}t\,,q_1(t)=-\frac{70}{243}t^9 + \frac{97}{162}t^6 - \frac{631}{1944}t^3 + \frac{1}{72}\,,\dots$. After specialization to~$t=0$, this gives again the Laurent expansion~\eqref{eq:ex1}.

\end{example}

\section{Laurent coefficients of regularized theta lifts.}
\label{sec:4}
In this section, we will prove Theorems \ref{thm:main} and \ref{thm:main2}, from which Theorem \ref{thm:intro} follows. Then we give an example in the Shimura curve case. 

\subsection{Proofs of Theorem \ref{thm:intro}}
For $f \in M^!_{k + \frac{1}{2}, L^-}$ with $k \ge 1$, let $\Phi_\BS(f, z, h)$ be the associated Borcherds-Shimura lift, which is meromorphic modular form of weight $2k$.
Let $a_m$ be the $m$-th Laurent coefficient of $\Phi_\BS$ at  a CM point $[z_U^\pm, t_0] \in Z(U) \subset X_K$, and $\tilde\Theta_N^{(k+m)} \in H_{1-k-m, L^-}$
a harmonic Maass form with $\xi$-image $\Theta_{N}^{(k+m)}$. 
Then we have the following result.
\begin{theorem}
  \label{thm:main}
    Write $m = 2r + \ell \ge 0$ with $\ell \in\{ 0, 1\}$, and let $P = L \cap (U^-)^\perp, N = L \cap U^-$. 
  In the notations above, we have
\begin{equation}
  \label{eq:main}
  a_m =
(-1)^{r+1}  \frac{(4\pi \sqrt{2}y_U)^m}{m! \binom{-k -m + r}{r}}
   \CT( \langle f(\tau), [\tilde\Theta^{k+m, +}_{N}(\tau, t_0), (2\sqrt{\pi})^{-\ell}\Theta^{\ell}_{P}(\tau)]_r\rangle).
\end{equation}
\end{theorem}
\begin{remark}
  \label{rmk:simplify}
  As $\ell \in \{0, 1\}$, the theta series $\Theta_P^{(\ell)}$ is holomorphic of weight $\ell + \frac12$.
  Using \eqref{eq:CT-RC}, one can rewrite
  $$
  \CT( \langle f, [\tilde\Theta^{(k+m), +}_{N}, \Theta^{(\ell)}_{P}]_r\rangle)
= \CT( \langle [\Theta^{(\ell)}_P, f]_r, \tilde\Theta^{(k+m), +}_{N}\rangle),
$$
which can be expressed as the regularized inner product between the weakly holomorphic modular form $[\Theta^{(\ell)}_P, f]_r$ and the holomorphic binary theta series $\Theta^{(k+m)}$.
If we denote $c^{(\ell)}_{P+\mu}(n)$ and $c^{(k'),+}_{N+\nu}(n)$ the $n$-th Fourier coefficient of $\langle \Theta_P^{(\ell)}, \phi_\mu\rangle$ and
$\langle \tilde\Theta_N^{(k'), +}, \phi_\nu\rangle$ respectively, then the constant term above can be explicitly written as
\begin{equation}
  \label{eq:CT-explicit}
  \begin{split}
&      \CT( \langle f, [\tilde\Theta^{(k+m), +}_{N}, \Theta^{(\ell)}_{P}]_r\rangle)\\
&  = \sum_{\substack{n_1, n_2 \in \Qb\\ (\mu, \nu) \in N'/N  \oplus P'/P}}
    c^{(k+m), +}_{N+\mu}(n_1)
    c^{(\ell)}_{P+\nu}(n_2)
    \varphi_{f, -n_1 - n_2}(\mu, \nu)
    (-n_1-n_2)^r P_r^{(\ell - 1/2, -k-m)} \lp \frac{n_1 - n_2}{n_1 + n_2} \rp
  \end{split}
\end{equation}
using \eqref{eq:RCq} with
$P^{(\alpha, \beta)}_r$ the Jacobi polynomial and
$\varphi_{f, m} \in \Sc_L \subset \Sc(\hat V)$ defined in \eqref{eq:varphim}. 
\end{remark}

\begin{proof}
  We apply Propositions \ref{prop:LaurentExp}, \ref{prop:over-reg} and \ref{prop:diff-id}
  to obtain
  \begin{align*}
  \frac{  m!a_m}{y_U^m} &=    \CT^{(2k+2m)}_{z_U} (\Rz^m \Phi_\BS(f, z, h))
          =
(2\sqrt{2\pi})^m
          \CT^{(2k+2m)}_{z_U} ( \Phi^{(m, k+m)}(f, z,h)) \\
&          =
(2\sqrt{2\pi})^m
           \Phi^{(m, k+m)}(f, z_U,t_0) 
        =
(2\sqrt{2\pi})^m(-1)^r
          \int^\reg \langle f(\tau), R_\tau^{r} \Theta^{(\ell, k+m)}_L(\tau, (z_U, t_0)) \rangle
  \end{align*}
  Using equations \eqref{eq:theta-split}, \eqref{eq:sub-pair} and applying the Rankin-Cohen bracket along with \eqref{eq:RCid}, we can write
  \begin{align*}
\langle f(\tau), R_\tau^{r} \Theta^{(\ell, k+m)}_L(\tau, (z_U, t_0)) \rangle   & = \left\langle f(\tau), R_\tau^{r}
           \lp \Theta^{(\ell)}_{P}(\tau)
                                                                                  \otimes \overline{\Theta^{(k+m)}_{N}(\tau, t_0)} v^{k+m+1}\rp\right\rangle\\
                 &= \left\langle f(\tau), R_\tau^r L_\tau 
           \lp \Theta^{(\ell)}_{P}(\tau)
           \otimes \tilde\Theta^{(k+m)}_{N}(\tau, t_0)\rp \right\rangle\\
                   &=
             (4\pi)^{r} \binom{-k - r - \ell}{r}^{-1}
                   \langle f(\tau),
             L_\tau [ \tilde\Theta^{(k+m)}_{N}(\tau, t_0), \Theta^{(\ell)}_{P}(\tau)]_r \rangle    .
  \end{align*}
  Applying Stokes' theorem (see e.g.\ \cite[Lemma 2.3]{LZ22}) gives us
  \begin{align*}
    \int^\reg
                   \langle f(\tau),
             L_\tau [ \tilde\Theta^{(k+m)}_{N}(\tau, t_0), \Theta^{(\ell)}_{P}(\tau)]_r \rangle    
    &=
-      \CT( \langle f(\tau), [\tilde\Theta^{(k+m), +}_{N}(\tau, t_0), \Theta^{(\ell)}_{P}(\tau)]_r\rangle).
  \end{align*}
  Putting these together proves \eqref{eq:main}.
\end{proof}

\begin{proof}[Proof of Theorem \ref{thm:intro}]
  We specialize Theorem \ref{thm:main} to this case by taking $L = M_2(\Qb)^0 \cap M_2(\Zb)$
with the quadratic form $Q(X) = \det(X)$.
For a CM point $z_0 = \frac{-B + \sqrt{-D}}{2A}$ with  fundamental discriminant $-D < 0$, 
we have the sublattices
$$
P = \Zb \cdot \smat{-B}{-2C}{2A}{B},~
N^- = \Zb \cdot \smat{A}{B}{0}{-A} + \Zb \cdot \smat{0}{-C}{-A}{0},
$$
and we have the rational splitting $P \oplus N^- \subset L$.
The map $x\cdot \smat{A}{B}{0}{-A} + y\cdot \smat{0}{-C}{-A}{0} \stackrel{\iota}{\mapsto} A(x + y z_0)$ gives an isometry between $N$ and $(\af, \Nm)$.
It is straightforward to verify that
\begin{align*}
  \sqrt{-D}(x + yz_0) &= x' + y'z_0,~ \pmat{B}{-2C}{2A}{-B} \binom{x}{y} = \binom{x'}{y'},\\
(X, Z(z_0)) &= -\frac{\sqrt{-D}}A \iota(X).
\end{align*}
Now we specialize to $A = 1$ and $D$ odd. Then it is easy to check that for $\ell = 0, 1$, 
$$
\Theta^{(\ell)}_P = \sum_{b \in \Zb} \phi_{\Zb + \frac{b}{2D}} (b \sqrt{\pi / D})^\ell q^{b^2/(4D)},~ 
\Theta_N^{(k)}(\tau) = \frac{2^{3k/2}}{D^k} \sum_{\mu \in \Zb/D\Zb} \phi_{\Oc + \frac{\mu}{\sqrt{-D}}} \theta_\mu^{(k)}(\tau)
$$
Now $(\mu, b) \in \Zb/D\Zb \oplus \Zb/2\Zb \cong (N^-)'/N^- \oplus P'/P$ is in $L'/(N^- \oplus P) $ if and only if
\begin{equation}
  \label{eq:cond1}
\frac{b}{2D} \cdot   \pmat{-B}{-2C}{2}{B} + \iota^{-1}\lp  \frac\mu{D} (-B - 2 z_0) \rp \in L'/(N^- \oplus P).  
\end{equation}
Straightforward calculations yield
$$
\frac{b}{2} \cdot   \pmat{-B}{-2C}{2}{B} + \mu \cdot \iota^{-1}\lp   (-B - 2 z_0) \rp
\equiv (\mu + b/2) \pmat{-B}{-C}{1}{B}
\bmod{D}.
$$
So \eqref{eq:cond1} holds if and only if $\mu \equiv -b/2 \bmod{D}$. 
Using the well-known isomorphism between $M^{!, +}_{k+1/2}(4)$ and $M^{!}_{k+1/2, L^-}$, we can view $f$ as a vector-valued modular form. Applying Theorem \ref{thm:main}, Remark \ref{rmk:simplify} and setting $n_1 = \frac{n}D, n_2 = \frac{b^2}{4D}$ gives us
\begin{align*}
  a_m &=
  (-1)^{r+1}  \frac{(2\pi \sqrt{2}\sqrt{D})^m}{m! \binom{-k -m + r}{r}} (2\sqrt{\pi})^{-\ell}
   \CT( \langle f(\tau), [\tilde\Theta^{k+m, +}_{N}(\tau, t_0), \Theta^{\ell}_{P}(\tau)]_r\rangle)\\
      &=
-  \frac{(2\pi \sqrt{2D})^m}{m! \binom{-k -m + r}{r}} (2\sqrt{\pi})^{-\ell}
        \sum_{n, b \in \Zb}
  \frac{2^{3(k+m)/2}}{D^{k+m}}
  c^{(k+m), +}_{-b/2}\lp \frac{n}{D}\rp\\
&\times (b\sqrt{\pi/D})^\ell
c_f\lp \frac{-4n-b^2}{D}\rp
    (4D)^{-r}(4n + b^2)^r P_r^{(\ell - 1/2, -k-m)} \lp \frac{4n - b^2}{4n + b^2} \rp,
\end{align*}
which is equivalent to \eqref{eq:am-intro}.
\end{proof}

For $f \in M^!_{-k+3/2, L^-}$ with $k \ge 1$, we can give a similar expression for the Laurent coefficients of the Borcherds-Zemel lift $\Phi_\BZ(f, z, h)$ at a CM point.

\begin{theorem}
  \label{thm:main2}
  Let $m \ge 0$ and denote $a_m$ the $m$-th Laurent coefficient of $\Phi_\BZ(f, z, h)$ at $[z_U, t_0] \in Z(U)$.
  Then we have
  \begin{equation}
    \label{eq:main2}
    a_m =
      \frac{(2\sqrt{2\pi}y_U)^m(-1)^{r'+1}             (4\pi)^{r'} }{m! \binom{k+\ell-2-r'}{r'}}
       \CT( \langle f(\tau), [\tilde\Theta^{(k+m), +}_{N}(\tau, t_0), \Theta^{(1-\ell)}_{P}(\tau)]_{r'}\rangle).
  \end{equation}
where  $\ell, r, P, N$ are the same as in Theorem \ref{thm:main} and $r' := k + \ell - 1 + r$.
\end{theorem}

\begin{proof}
  As in the proof of Theorem \ref{thm:main}, we have
  \begin{align*}
m!    \frac{a_m}{y_U^m}
    &=
(2\sqrt{2\pi})^m
           \Phi^{(2k-1+m, k+m)}(f, z_U,t_0) \\
        &=
(2\sqrt{2\pi})^m(-1)^{k+r+\ell-1}
          \int^\reg \langle f(\tau), R_\tau^{k+r+\ell-1} \Theta^{(1-\ell, k+m)}_L(\tau, (z_U, t_0)) \rangle\\
    &=
      \frac{(2\sqrt{2\pi})^m(-1)^{k+r+\ell}             (4\pi)^{k+r+\ell-1} }{\binom{-r-1}{k + r + \ell-1}}
       \CT( \langle f(\tau), [\tilde\Theta^{(k+m), +}_{N}(\tau, t_0), \Theta^{(1-\ell)}_{P}(\tau)]_{k+r+\ell-1}\rangle).
  \end{align*}
  This finishes the proof.
\end{proof}

\begin{remark}
  Theorem \ref{thm:main2} can also be specialized to a more explicit form as in Remark \ref{rmk:simplify} and in Theorem \ref{thm:intro}. If we take the lattice as in the proof of Theorem \ref{thm:intro}, then the Borcherds-Zemel lift has the following Fourier expansion
\begin{equation}
  \label{eq:BZ-FE}
  \Phi_\BZ(f, z) =
    2^{1+k/2} \pi^{k-1/2} i^{1-k}
\lp C +   \sum_{n \ge 1} \lp \sum_{d \mid n} d^{k} c_f(d^2) \rp n^{k-1} q^n \rp,
\end{equation}
where $C$ is the constant term.

\end{remark}


\subsection{Example: Modular forms on Shimura curves}
\label{subsec:ex2}
Let $B$ be a quaternion algebra over $\Qb$ with discriminant $6$ \cite[\textsection 3]{Voight09}. It is the $\Qb$-algebra generated by $\alpha, \beta$ satisfying
\begin{equation}
  \label{eq:ab}
  \alpha^2 = 3,  \beta^2 = -1,  \gamma := \alpha \beta = - \beta \alpha.
\end{equation}
It is more convenient to view $B$ as a $\Qb$-subalgebra of $M_2(\Rb)$ via the embedding \footnote{We slightly abuse the notation by using $X$ to denote $\iota_\infty(X)$ for $X \in B_\Rb$.}
\begin{equation}
  \label{eq:iota_inf}
  \begin{split}
    \iota_\infty: B &\hookrightarrow M_2(\Rb) \\
    \alpha, \beta, \gamma &\mapsto  \sqrt{3} \pmat{1}{0}{0}{-1}, \pmat{0}{-1}{1}{0}
                     ,  -\sqrt{3} \pmat{0}{1}{1}{0}.
    .
  \end{split}
\end{equation}
The reduced norm $Q$ on $B$ then agrees with $\det$ on $M_2(\Rb)$.
The subspace $B^0 \subset B$ of trace zero elements is given by $\Qb \alpha \oplus \Qb \beta \oplus \Qb \gamma$, which is acted on by elements in $B$ with reduced norm 1 through conjugation.
Denote $\Oc := \Zb \oplus \Zb \alpha \oplus \Zb \beta \oplus \Zb \delta$  a maximal order with units $\Oc^*_1 \subset \Oc$. 
The Shimura curve $X^B(1) := \Oc_1^* \backslash \Hb$ has genus zero.

Following \cite{BG12}, the ring of holomorphic modular forms on $X^B(1)$ is generated by $\Pc, \Qc, \Rc$ of weight $4, 6, 12$ respectively.
They satisfy the relation $\Pc^6 + 3\Qc^4 + \Rc^2 = 0$.
In \cite[section 5]{BG12}, the Laurent expansion of $\Pc(z)$ at $z = i$ is given as
$$
(1-w)^{-4}\Pc(z) = i \lp \sqrt{2\pi} \Omega \rp^4
\lp \frac1{2^33^2} + \frac5{2^43^3} \frac{u^4}{4!} - \frac{5 \cdot 17}{2^4 3^4} \frac{u^8}{8!} + \frac{5 \cdot 317}{2^53^3} \frac{u^{12}}{12!} + \dots \rp,
$$
where $w = \frac{z-i}{z+i}$, $\Omega = \Omega_{-4} = \frac1{2\sqrt{2\pi}} \Gamma(\frac14)\Gamma(\frac34)^{-1} = 0.590170...$, 
and $u = 4\pi \Omega^2 w$. 
The goal is to obtain $\Pc$ as the Borcherds-Shimura lift of a vector-valued cusp form of weight $5/2$. 

For this, we take as our lattice
\begin{equation}
  \label{eq:L_compact}
L := B^0 \cap \Oc = \Zb \alpha \oplus \Zb \gamma \oplus \Zb \beta.
\end{equation}
Then $Q(\beta) = 1, Q(\alpha) = Q(\gamma) = -3$ and 
\begin{equation}
\label{eq:Ldual}
L' = \frac{1}{6} \Zb \alpha \oplus \frac{1}{6} \Zb \gamma \oplus  \frac{1}{2} \Zb \beta.
\end{equation}
The group $\Oc^*_1 \subset \SO(L)$ 
is generated by
$$
-1, h_3 := \beta, h_2 := (-1 - \alpha + 3\beta +\gamma)/2,
h_2' := (-1 + \alpha + 3\beta +\gamma)/2.
$$
The conjugation action of $\Oc^*_1$ in the basis $\{\alpha, \beta, \gamma\}$ is given by
$$
h_3 =
\begin{pmatrix}
  -1 & & \\
  & 1 & \\
  & & -1
\end{pmatrix},~
h_2 =
\begin{pmatrix}
  -2 & -1 & 0 \\
6 & 4 & -3 \\
3 & 2 & -2
\end{pmatrix},~
h_2' =
\begin{pmatrix}
  -2 & 2 & -3 \\
-3 & 4 & -6 \\
0 & 1 & -2
\end{pmatrix}.
$$
This fixes $L$ and its image $H_1$ in the group $H_0 := \SO(L'/L)$ is generated by
\begin{equation}
  \label{eq:sigma}
  h_1 := h_2 h_3
\end{equation}
and isomorphic to $\Zb/6\Zb$. 
To produce modular forms that lifts to $X^B(1)$, we need to lift from the $H_1$-invariant subspace.
Thus, we first need to analyze and decompose $\rho^{H_1}_{L^-}$.

We write $A:= L'/L = \{(a, c, b): a, c \in \Zb/6\Zb, b \in \Zb/2\Zb\}$ with quadratic form
$$
Q(a, c, b) = \frac{-a^2 - c^2+ 3b^2 }{12}
$$
%
%
Denote $A_p := A \otimes \Zb_p$ for $p = 2, 3$.
Then $A_2 = (\Zb/2\Zb)^3$ and $A_3 = (\Zb/3\Zb)^2$ with quadratic forms
$$
Q_2(a_2, c_2, b_2) = \frac{a_2^2 + c_2^2 + b_2^2}{4} ,~
Q_3(a_3, c_3) = \frac{-a_3^2 - c_3^2}{3}
$$
valued in $\Qb/\Zb$.
We identify
\begin{equation}
  \label{eq:A23}
 A \cong A_2 \oplus A_3,~ (a, c, b) \mapsto ((a, c, b) \bmod{2} , (a, c) \bmod{3}).
\end{equation}
Using \eqref{eq:A23}, we can identify $\SO(A_p)$ as the subgroup of $H_0 = \SO(A)$, from which we obtain $\SO(A) \cong \SO(A_2) \times \SO(A_3)$.
Note that the group $\SO(A_p)$ is isomorphic to the dihedral group of order $2p+2$, and generated by $\{\sigma_p, \tau_p\}$ with
\begin{equation}
  \label{eq:sigma-tau}
  \begin{split}
    \sigma_2(a_2, c_2, b_2) &= (c_2, b_2, a_2),~
  \tau_2(a_2, c_2, b_2) = (c_2, a_2, b_2),\\
  \sigma_3(a_3, c_3) &= (c_3, -a_3),~
  \tau_3(a_3, c_3) = (c_3, a_3).
  \end{split}
\end{equation}
For $p = 2, 3$, the element $h_p$ is given by
$$
h_3 = \sigma_3^2,~
h_2 = \sigma_2,~
h_2' = \sigma_2^{-1}.
$$
For $\mu = (\mu_2, \mu_3) \in A$ with $\mu_p \in A_p$, it is easy to see that
$$
\Oc_1^* \cdot \mu =
\{(\sigma_2^i \mu_2, \sigma_3^{2j} \mu_3): i = 0, 1, 2, j = 0, 1\}. 
$$
We write $H_1 = H_2 \times H_3$ with $H_p = \langle h_p \rangle$ for $p = 2, 3$.
For $l = 1, 2$, we define characters $\chi_l$ of $H_0$ by
\begin{equation}
  \label{eq:chij}
  \begin{split}
    \chi_{1}(\sigma_2) &= 1,~
      \chi_{1}(\tau_2) = 1,~ 
    \chi_{1}(\sigma_3) = -1,~
      \chi_{1}(\tau_3) = -1,\\
    \chi_{2}(\sigma_2) &= 1,~
      \chi_{2}(\tau_2) = 1,~ 
                         \chi_{2}(\sigma_3) = -1,~  \chi_{2}(\tau_3) = 1.
  \end{split}
\end{equation}
Then it is straightforward to verify that
\begin{equation}
  \label{eq:CAdecomp}
  \Cb[A]^{H_1} = \Cb[A]^{H_0} \oplus \Cb[A]^{\chi_1} \oplus \Cb[A]^{\chi_2}.
\end{equation}
with $\Cb[A]^{H_i}$, resp.\ $\Cb[A]^{\chi_l}$, the $H_i$-invariant, resp.\ $\chi_l$ isotypic, subspace of $\Cb[A]$.
      For $l =1, 2$, we also have the projection operator $\Pi_l$ to the $\chi_l$-isotypic component defined by \footnote{We also denote $\Pi_0$ the projection to $H_0$-invariant subspace using the trivial character of $H_0$.}
      \begin{equation}
        \label{eq:Pil}
        \begin{split}
                  \Pi_l(\phi) &:= \frac1{|\SO(A)|} \sum_{h \in \SO(A)} \chi_l(h)^{-1} (h\cdot \phi).
        \end{split}
      \end{equation}
Extending this decomposition to functions valued in $\Cb[A]$ gives us
\begin{equation}
\label{eq:vvmf-isom}
M^!_{k+1/2, \rho} =
M^{!, H_1}_{k+1/2, L^-} =
	M^{!, H_0}_{k+1/2, L^-} \oplus
	M^{!, \chi_1}_{k+1/2, L^- } \oplus
	M^{!, \chi_2}_{k+1/2, L^-}.
      \end{equation}
      where 
      the superscript denotes the $H_i$-invariant or $\chi_l$-isotyptic subspace.

%

For any weakly holomorphic modular form $f$ in $M^!_{k+1/2, \rho}$, resp.\ $M^!_{3/2-k, \rho}$, the modular form $\Phi_\BS(f, z)$, resp.\ $\Phi_\BZ(f, z)$, is a meromorphic modular form of weight $2k$ on $X^B(1)$. 
To find the cusp form that lifts to $\Pc$, we need to look at the space 
$S^{H_0}_{5/2, L^-}$ and $S^{\chi_l}_{5/2, L^-}$ for $l = 1, 2$. 
Using the SAGE code WeilRep by Brandon Williams \cite{WeilRep}, we quickly see that these 3 spaces have dimensions 1, 1, 0 respectively. 
To see which form $f$ lifts to a non-trivial multiple of $\Pc$, we can check the value of $\Phi_\BS(f, z)$
%
at $z = i$.

At the point $z = i$, we have $U = \Qb \alpha + \Qb \gamma$,
$N^- = U \cap L = \Zb \alpha + \Zb \gamma$,  $P = \Zb \beta$ and $A = B \oplus C$ with $B := P'/P, C := N'/N$.
This gives $\Cb[A] = \Cb[B] \otimes \Cb[C]$.
      By viewing $\chi_l$ as a character on the subgroup  $\SO(C) \subset \SO(A)$, which is generated by $\tau_2, \sigma_3, \tau_3$, we can define a projection operator $\Pi'_l$ on $\Cb[C]$ as in \eqref{eq:Pil} by
      \begin{equation}
        \label{eq:Pil'}
        \begin{split}
                  \Pi'_l(\phi) &:= \frac1{|\SO(C)|} \sum_{h \in \SO(C)} \chi_l(h)^{-1} (h\cdot \phi).
        \end{split}
      \end{equation}
      Also, it is clear that $        \Pi'_l(\phi) = \phi$
      for any $\phi \in \Cb[A]^{\chi_l}$.
      In particular, we have 
      \begin{equation}
        \label{eq:Pil'-equiv}
        \langle f, \Theta_L \rangle
= \langle \Pi'_l(f), \Theta_L \rangle         
        = \langle f, \Pi'_l(\Theta_L) \rangle 
\end{equation}
for $f \in M^{!, \chi_l}_{k+1/2, L^-}$.

The theta function splits as
$$
\Theta^{(\ell, k)}_L(\tau, i) = \Theta^{(\ell)}_P(\tau) \otimes
\overline{\Theta^{(k)}_{N}(\tau) } v^{k+1}
$$
with
\begin{align*}
  \Theta^{(\ell)}_P(\tau) 
  &= (4v)^{-\ell/2} \sum_{b \in \Zb/2\Zb} \phi_{b/2} \sum_{\lambda \in P + \frac{b}2 \beta}
    \He_\ell(\sqrt{\pi v} (\lambda, \beta)) q^{Q(\lambda)},\\
  \Theta^{(k)}_{N}(\tau)
  &= (\sqrt6 i)^{k} \sum_{\nu \in N'/N} \phi_\nu \sum_{a\alpha + c \gamma \in N + \nu}
    (a + ci)^k q^{3(a^2 + c^2)} \in S_{k+1, N}.
\end{align*}
To obtain the value of $\Phi_\BS(f, i)$, we need to know $\Pi'_l(\Theta^{(0, 2)}_L(\tau, i))$ for $l = 0, 1, 2$ by \eqref{eq:Pil'-equiv}.
A quick calculation using SAGE shows that  $\Pi'_0(\Theta_N^{(2)})$ 
vanishes identically.
Therefore $\Phi_\BS(f, i) = 0$ for $f \in S^{H_0}_{5/2, L^-}$ and $\Phi_\BS(f, z)$ is identically 0 for such $f$.
Since $S^{\chi_2}_{5/2, L^-}$ is trivial, the only forms that could lift to $\Pc$ lie in $S^{\chi_1}_{5/2, L^-}$.
So we focus on the $\chi_1$-isotypic part of $\Cb[A] \cong \Cb[A_2] \otimes \Cb[A_3]$, which has the following basis $\{\ef_i: 0 \le i \le 3\}$ with
\begin{equation}
  \label{eq:efs}
  \begin{split}
    \ef_i &:= \sum_{a, b, c \in \{0, 1\},~ a + b + c = i} \phi_{(a/2, c/2, b/2)} \otimes (\phi_{(0, 1)} + \phi_{(0, 2)} - \phi_{( 1,0)}  - \phi_{( 2,0)}) \in \Cb[A]^{\chi_1}.
  \end{split}
\end{equation}
Later, we will also look at the $\chi_1$-isotypic part of $\Cb[C]$, which has the basis $\{\ef_m': 0 \le m \le 2\}$ with
\begin{equation}
  \label{eq:efs}
  \begin{split}
    \ef'_m &:= \sum_{a, c \in \{0, 1\},~ a + c = m} \phi_{(a/2, c/2)} \otimes (\phi_{(0, 1)} + \phi_{(0, 2)} - \phi_{( 1,0)}  - \phi_{( 2,0)}) \in \Cb[C]^{\chi_1}.
  \end{split}
\end{equation}
Note that the space $S^{\chi_1}_{5/2, L^-}$ is spanned by
\begin{equation}
  \label{eq:f1}
 f_1 = q^{1/12} \ef_1 + 6 q^{1/3} \ef_0 + O(q^{5/6})\in S^{\chi_1}_{5/2, L^-}.
\end{equation}

Consider $\Pi_1'(\Theta_N^{(k)})$ in the space $S_{k+1, N}^{\chi_1}$, which is trivial when $k \neq 2 \bmod4$. 
We can express them explicitly using Poincar\'e series via
\begin{equation}
  \label{eq:Pseries}
  \begin{split}
    \Pi'_1(\Theta^{(2)}_N)
    &= \sqrt{3} \pi^2\Omega^4  \frac{1}{3^3} P_{3, 1/12, \ef'_1}
= \frac1{12} (q^{1/12} \ef'_1 + 8 q^{1/3} \ef'_0 - 32 q^{5/6}\ef'_2 + O(q^{13/12}))      ,\\
    \Pi'_1(\Theta^{(6)}_N)
    &= \sqrt3 \pi^6 \Omega^{12}  \frac{ 2^2}{3^9 \cdot 5} (19 P_{7, 1/12, \ef'_1} + 2^{11} P_{7, 1/3, \ef'_0})\\
&= \frac{1}{2^4 \cdot 3^3} (q^{1/12} \ef'_1 + 128q^{1/3} \ef'_0   + 1408 q^{5/6}\ef'_2 +O(q^{13/12}))      .
  \end{split}
\end{equation}
Here $P_{k, m, \phi}$ is the usual vector-valued Poincar\'e series with $\phi \in \Cb[C]$ (see \cite{bruinierhabil}).
Let $F_{k, m, \phi}$ be the Maass-Poincar\'e series \cite{bruinierhabil}, which is a harmonic Maass form satisfying
$$
\xi_k F_{k, m, \phi} = \frac{(4\pi |m|)^{1-k}}{(-k)!} P_{2-k, -m, \phi},~
 F^+_{k, m, \phi} = q^{m} \phi + O(1),
 $$
 where $F^+_{k, m, \phi}$ is the holomorphic part of $F_{k, m, \phi}$.
 Therefore, \eqref{eq:Pseries} tells us that
 \begin{equation}
   \label{eq:tThetaN}
   \begin{split}
     \tilde\Theta_N^{(2)}
     &:= \frac{\sqrt3 }3 \Omega^4 F_{-1, -1/12, \ef'_1},~
     \tilde\Theta_N^{(6)}
     := \frac{2^4\sqrt3 }{3^2} \Omega^{12} (38 F_{-5, -1/12, \ef'_1} + F_{-5, -1/3, \ef'_0})
   \end{split}
 \end{equation}
 are $\xi$-preimages of $\Pi'_1(\Theta_N^{(k)})$ for $k = 2, 6$.
 Their holomorphic parts are given by
 \begin{equation}
   \label{eq:tThetaN+}
   \begin{split}
     \tilde\Theta_N^{(2), +}
     &= \frac{\sqrt3 }3 \Omega^4 (q^{-1/12} \ef'_1 + O(1)),~ 
     \tilde\Theta_N^{(6), +}
       = \frac{2^4\sqrt3 }{3^2} \Omega^{12}(q^{-1/3} \ef'_0 + 38 q^{-1/12} \ef'_1+ O(1)).
   \end{split}
 \end{equation}

By Theorem \ref{thm:main}, we then have
\begin{align*}
&  \Phi_\BS(f_1, i)
  = a_0 = - \CT \langle f_1, \tilde\Theta^{(2), +}_N \otimes \Theta_P\rangle\\
  &  =      - \CT \sum_{i, b, m}
    \frac{ \langle f_{1}, \ef_i\rangle }{\langle \ef_i, \ef_i\rangle}
\frac{   \langle \tilde\Theta^{(2), +}_{N}, \ef'_m\rangle}{\langle \ef'_m, \ef'_m\rangle}
    \langle     \Theta_P ,  \phi_{b/2} \rangle
\langle \ef_i,    \ef'_m \otimes \phi_{b/2}\rangle\\
&= - \CT
    \frac{ \langle f_{1}, \ef_1\rangle }{\langle \ef_1, \ef_1\rangle}
\frac{   \langle \tilde\Theta^{(2), +}_{N}, \ef'_1\rangle}{\langle \ef'_1, \ef'_m\rangle}
    \langle     \Theta_P ,  \phi_{b/2} \rangle
\langle \ef_1,    \ef'_1 \otimes \phi_1\rangle
       = - \frac{8 \sqrt3 \Omega^4 }{3}. 
    \end{align*}
    Therefore, we have $\Phi_\BS(f_1, z) = \frac{2\sqrt3 i }{3\pi^2} \Pc(z)$. 
    As a check, we can apply Theorem \ref{thm:main} and \eqref{eq:Pil'-equiv} to find the 4-th Laurent coefficient $a_4$ of $\Phi_{\BS}(f_1, z)$ at $z = i$ is given by
\begin{align*}
  a_4 &=  (-1)^{5}  \frac{(4\pi \sqrt{2})^4}{4! \binom{-2 -4 + 2}{2}}\CT \langle f_1, ( [\tilde\Theta^{(6), +}_N , \Theta_P]_2)\rangle\\
      &=      - \frac{2^6 \pi^4}{15} \CT \sum_{i, b, m, s}
 (-1)^s \binom{-4}s \binom{3/2}{2-s}
    \frac{ \langle f_{1}, \ef_i\rangle }{\langle \ef_i, \ef_i\rangle}
        \frac{   \langle (q\partial_q)^{2-s} \tilde\Theta^{(2), +}_{N}, \ef'_m\rangle}{\langle \ef'_m, \ef'_m\rangle}
    \langle (q\partial_q)^{s}     \Theta_P ,  \phi_{b/2} \rangle
        \langle     \ef_i ,  \ef'_m \otimes \phi_{b/2} \rangle\\
      &=      - \frac{2^6 \pi^4}{15} \CT \sum_{(i, b, m) \in S_0}
\frac38
    \frac{ \langle f_{1}, \ef_i\rangle }{\langle \ef_i, \ef_i\rangle}
        \frac{   \langle (q\partial_q)^{2} \tilde\Theta^{(2), +}_{N}, \ef'_m\rangle}{\langle \ef'_m, \ef'_m\rangle}
    \langle      \Theta_P ,  \phi_{b/2} \rangle
        \langle     \ef_i ,  \ef'_m \otimes \phi_{b/2} \rangle\\
      & \quad
     - \frac{2^6 \pi^4}{15} \CT \sum_{s = 1, 2}
        (-1)^s \binom{-4}s \binom{3/2}{2-s}
    \frac{ \langle f_{1}, \ef_1\rangle }{\langle \ef_1, \ef_1\rangle}
        \frac{   \langle (q\partial_q)^{2-s} \tilde\Theta^{(2), +}_{N}, \ef'_0\rangle}{\langle \ef'_0, \ef'_0\rangle}
    \langle (q\partial_q)^{s}     \Theta_P ,  \phi_1 \rangle
        \langle     \ef_1 ,  \ef'_0 \otimes \phi_1 \rangle\\
      &=        - \frac{2^{10}\sqrt3  \pi^4 \Omega^{12}}{3^3 \cdot 5}
        \lp \frac38 \cdot \lp 6 \cdot \frac19 \cdot 4 + 38\cdot \frac1{144} \cdot 8 + \frac29 \cdot 4 \rp
+ 6 \cdot \frac{-1}6 \cdot 4 + 10\cdot \frac1{8} \cdot 4
\rp
        =       - \frac{2^7 \cdot 5 \sqrt3  \pi^4 \Omega^{12}}{3^3}.
\end{align*}
Here $S_0 = \{(0, 0, 0), (1, 1, 0), (1, 0, 1)\}$.
Then we see that
$$
\lp \frac{2i}{z+i}\rp^{-4}
\Phi_\BS(f_1, z) = a_0 + a_4w^4 + O(w^8) =
 - \frac{8 \sqrt3 \Omega^4 }{3}\lp 1 + \frac56 \frac{ (4 \pi \Omega^2 w)^4}{4!} + O(w^8)\rp,
 $$
 which agrees with the Laurent expansion of $\frac{2\sqrt3 i }{3\pi^2} \Pc(z)$ up to $O(w^8)$.

\bibliography{TECM}{}

\def\cprime{$'$}
\begin{thebibliography}{BvdGHZ08}

\bibitem[BF04]{bruinierfunke04}
Jan~H. Bruinier and Jens Funke.
\newblock On two geometric theta lifts.
\newblock {\em Duke Math. J.}, 125(1):45--90, 2004.

\bibitem[BG12]{BG12}
Srinath Baba and H{\aa}kan Granath.
\newblock Differential equations and expansions for quaternionic modular forms
  in the discriminant 6 case.
\newblock {\em LMS J. Comput. Math.}, 15:385--399, 2012.

\bibitem[BKY12]{BKY12}
Jan~Hendrik Bruinier, Stephen~S. Kudla, and Tonghai Yang.
\newblock Special values of {G}reen functions at big {CM} points.
\newblock {\em Int. Math. Res. Not. IMRN}, (9):1917--1967, 2012.

\bibitem[Bor98]{borcherds}
Richard~E. Borcherds.
\newblock Automorphic forms with singularities on {G}rassmannians.
\newblock {\em Invent. Math.}, 132(3):491--562, 1998.

\bibitem[Bru02]{bruinierhabil}
Jan~H. Bruinier.
\newblock {\em {B}orcherds products on $O(2,l)$ and {C}hern classes of
  {H}eegner divisors}.
\newblock Springer Lecture Notes in Mathematics 1780, 2002.

\bibitem[BvdGHZ08]{brdgza08}
Jan~H. Bruinier, Gerard van~der Geer, G{\"u}nter Harder, and Don Zagier.
\newblock {\em The 1-2-3 of modular forms}.
\newblock Universitext. Springer-Verlag, Berlin, 2008.
\newblock Lectures from the Summer School on Modular Forms and their
  Applications held in Nordfjordeid, June 2004, Edited by Kristian Ranestad.

\bibitem[BY09]{bruinieryang}
Jan~H. Bruinier and Tonghai Yang.
\newblock {F}altings heights of {CM} cycles and derivatives of $l$-functions.
\newblock {\em Invent. Math.}, 177:631--681, 2009.

\bibitem[ELS22]{ehlenlischwagenscheidt}
Stephan Ehlen, Yingkun Li, and Markus Schwagenscheidt.
\newblock {H}armonic {M}aass forms associated with {CM} newforms.
\newblock {\em preprint arXiv:2210.07341}, 2022.

\bibitem[FM06]{FM06}
Jens Funke and John Millson.
\newblock Cycles with local coefficients for orthogonal groups and
  vector-valued {S}iegel modular forms.
\newblock {\em Amer. J. Math.}, 128(4):899--948, 2006.

\bibitem[GMR20]{guerzhoymertensrolen}
Pavel Guerzhoy, Michael Mertens, and Larry Rolen.
\newblock {P}eriodicity of {T}aylor coefficients of half-integer weight modular
  forms.
\newblock {\em Pacific Journal of Mathematics}, 307:137--157, 2020.

\bibitem[Li22]{Li22-average}
Yingkun Li.
\newblock Average {CM}-values of higher {G}reen's function and factorization.
\newblock {\em Amer. J. Math.}, 144(5):1241--1298, 2022.

\bibitem[LZ22]{LZ22}
Yingkun Li and Shaul Zemel.
\newblock Shintani lifts of nearly holomorphic modular forms.
\newblock {\em Canadian Journal of Mathematics}, pages 1--47, 2022.

\bibitem[Nel11]{nelson1}
Paul Nelson.
\newblock Computing on {S}himura curves.
\newblock Appendices B and C of Periods and special values of $L$-functions,
  notes from Arizona Winter School 2011, 2011.

\bibitem[Nel15]{nelson2}
Paul Nelson.
\newblock Evaluating modular forms on {S}himura curves.
\newblock {\em Mathematics of computation}, 84(295):2471--2503, 2015.

\bibitem[OR13]{osullivanrisager}
Cormac O'Sullivan and Morten~S. Risager.
\newblock Non-vanishing of {T}aylor coefficients and {P}oincar\'{e} series.
\newblock {\em Ramanujan J.}, 30(1):67--100, 2013.

\bibitem[Rom19]{romik}
Dan Romik.
\newblock {T}aylor expansion of the {J}acobi theta constant $\theta_3$.
\newblock {\em Ramanujan J.}, pages 1--16, 2019.

\bibitem[Sch09]{Schofer}
Jarad Schofer.
\newblock Borcherds forms and generalizations of singular moduli.
\newblock {\em J. Reine Angew. Math.}, 629:1--36, 2009.

\bibitem[Sch21]{scherer}
Robert Scherer.
\newblock Congruences modulo primes of the {R}omik sequence related to the
  {T}aylor expansion of the {J}acobi theta constant $\theta_3$.
\newblock {\em Ramanujan J.}, 54(427--448), 2021.

\bibitem[Voi09]{Voight09}
John Voight.
\newblock Shimura curve computations.
\newblock In {\em Arithmetic geometry}, volume~8 of {\em Clay Math. Proc.},
  pages 103--113. Amer. Math. Soc., Providence, RI, 2009.

\bibitem[VW14]{voightwillis}
John Voight and John Willis.
\newblock Computing power series expansions of modular forms.
\newblock In {\em Computations with modular forms}, volume~6 of {\em Contrib.
  Math. Comput. Sci.}, pages 331--361. Springer, Cham, 2014.

\bibitem[VZ93]{villegaszagier}
Fernando~Rodriguez Villegas and Don Zagier.
\newblock Square roots of central values of {H}ecke {$L$}-series.
\newblock In {\em Advances in number theory ({K}ingston, {ON}, 1991)}, Oxford
  Sci. Publ., pages 81--99. Oxford Univ. Press, New York, 1993.

\bibitem[Wil]{WeilRep}
Brandon Williams.
\newblock weilrep - {S}age code for computing with vector-valued modular forms,
  {J}acobi forms and theta lifts.
\newblock {\em \url{https://github.com/btw-47/weilrep}}.

\bibitem[Zem15]{Zemel15}
Shaul Zemel.
\newblock A {G}ross-{K}ohnen-{Z}agier type theorem for higher-codimensional
  {H}eegner cycles.
\newblock {\em Res. Number Theory}, 1:Paper No. 23, 44, 2015.

\end{thebibliography}
\bibliographystyle{alpha}

\end{document}